\documentclass[preprint,3p,12pt]{elsarticle}

\usepackage{framed,multirow}

\usepackage{amssymb}
\usepackage{latexsym}

\usepackage{url}
\usepackage{xcolor}
\definecolor{newcolor}{rgb}{.8,.349,.1}

\usepackage{bm}
\usepackage{tabu}
\usepackage{subfigure}
\usepackage{amsmath}
\usepackage{booktabs}

\DeclareMathOperator{\sign}{sign}

\begin{document}


\begin{frontmatter}

\title{PFNN: A Penalty-Free Neural Network Method for Solving a Class of Second-Order Boundary-Value Problems on Complex Geometries}

\author[iscas,ucas]{Hailong {Sheng}}
\ead{hailong2019@iscas.ac.cn}
\author[pku]{Chao {Yang}\corref{cor1}}
\cortext[cor1]{Corresponding author.}
\ead{chao_yang@pku.edu.cn}

\address[iscas]{Institute of Software, Chinese Academy of Sciences, Beijing 100190, China}
\address[pku]{School of Mathematical Sciences, Peking University, Beijing 100871, China}
\address[ucas]{University of Chinese Academy of Sciences, Beijing 100190, China}


\begin{abstract}
We present PFNN, a penalty-free neural network method, to efficiently solve  a class of second-order boundary-value problems on complex geometries. To reduce the smoothness requirement, the original problem is reformulated to a weak form so that the evaluations of high-order derivatives are avoided. Two neural networks, rather than just one, are employed to construct the approximate solution, with one network satisfying the essential boundary conditions and the other handling the rest part of the domain. In this way, an unconstrained optimization problem, instead of a constrained one, is solved without adding any penalty terms.
The entanglement of the two networks is eliminated with the help of a length factor function that is scale invariant and can adapt with complex geometries. We prove the convergence of the PFNN method and conduct numerical experiments on a series of linear and nonlinear  second-order boundary-value problems to demonstrate that PFNN is superior to several existing approaches in terms of accuracy, flexibility and robustness.
\end{abstract}

\begin{keyword}
  Deep neural network \sep
  Penalty-free method \sep
  Boundary-value problem \sep
  Partial differential equation \sep
  Complex geometry
\end{keyword}

\end{frontmatter}


\section{Introduction}
During the past decade, neural network is gaining increasingly more attention from the big data and artificial intelligence community \cite{goodfellow2016deep, lecun2015deep, schmidhuber2015deep}, and playing a central role in a broad range of applications related to, e.g., image processing \cite{razzak2018deep}, computer vision \cite{voulodimos2018deep} and natural language processing \cite{young2018recent}.
However, in scientific and engineering computing, the studying of neural networks is still at an early stage.
Recently, thanks to the introduction of deep neural networks, successes have been achieved in a series of challenging tasks, including turbulence modeling \cite{ling2016reynolds, duraisamy2019turbulence}, molecular dynamics simulations \cite{zhang2018deep, li2016understanding}, and solutions of stochastic and high-dimensional partial differential equations \cite{weinan2017deep, han2018solving, raissi2018forward, beck2019machine}.
In addition to that, a great deal of efforts are also made to build the connections between neural networks and traditional numerical methods such as finite element \cite{he2019relu}, wavelet \cite{fan2019bcr}, multigrid \cite{he2019mgnet}, hierarchical matrices \cite{fan2019multiscale1, fan2019multiscale2} and domain decomposition \cite{li2019d3m} methods.
However, it is still far from clear whether neural networks can solve ordinary/partial differential equations and truly exhibit advantages as compared to classical discretization schemes in accuracy, flexibility and robustness.

An early attempt on utilizing neural network methods to solve differential equations can be traced back to three decades ago \cite{lee1990neural}, in which a Hopfield neural network was employed to represent the discretized solution. Shortly afterwards, methodologies to construct closed form numerical solutions with neural networks  were proposed \cite{meade1994numerical, van1995neural}. Since then, more and more efforts were made in solving differential equations with various types of neural networks, such as feedforward neural network \cite{lagaris1998artificial, lagaris2000neural, mcfall2009artificial, mcfall2013automated}, radial basis network \cite{jianyu2002numerical, mai2001numerical}, finite element network \cite{ramuhalli2005finite}, cellular network \cite{chua1988cellular}, and wavelet network \cite{li2010integration}. These early studies have in certain ways illustrated the feasibility and effectiveness of neural network based methods, but are limited to handling model problems with regular geometry in low space dimensions. Moreover, the networks utilized in these works are relatively shallow, usually containing only one hidden layer, whereas the potential merits of the deep neural networks are not fully revealed.

Nowadays, with the advent of the deep learning technique \cite{nielsen2015neural, aggarwal2018neural, higham2019deep}, neural networks with substantially more hidden layers have become valuable assets. Various challenging problems modeled by complex differential equations have been taken into considerations and successfully handled by neural network based methods \cite{weinan2018deep, liao2019deep, zang2020weak, raissi2020hidden, raissi2019physics, sirignano2018dgm, berg2018unified, long2019pde, kani2017dr, khoo2017solving}. These remarkable progresses have  demonstrated the advantages of deep neural networks in terms of the strong representation capability. However, the effectiveness of the neural network based methods are still hindered by several factors. Specifically, the accuracy of the neural network based approximation as well as the efficiency of the training task usually have a strong dependency on the properties of the problems, such as the nonlinearity of the equation, the irregularity of the solution, the shape of the boundary, the dimension of the domain, among other factors. It is therefore of paramount importance to improve the robustness and flexibility of the neural network based approaches. For instance, several methods \cite{weinan2018deep, liao2019deep, zang2020weak, li2019d3m} have been recently proposed to transform the original problem into a corresponding weak form, thus lowering the smoothness requirement and possibly reducing the cost of the training process. Unfortunately, in these works extra penalty terms due to the boundary constraints are included in the loss function, which could eventually have a negative effect on the training process and the achievable accuracy. On the other hand, a number of efforts \cite{lagaris1998artificial, lagaris2000neural, mcfall2009artificial, mcfall2013automated} were made to explicitly manufacture the numerical solution so that the boundary conditions can be automatically satisfied. However, the construction processes are usually limited to problems with simple geometries, not flexible enough to generalize to arbitrary domains with arbitrary boundary conditions.

In this work, we propose PFNN, a penalty-free neural network method,
for solving a class of second-order boundary-value problems on complex geometries. 
In order to reduce the smoothness requirement, we reformulate the original problem to a weak form so that the approximation of high-order derivatives of the solution is avoided. In order to adapt with various boundary constrains without adding any penalty terms, we employ two networks, rather than just one, to construct the approximate solution, with one network satisfying the essential boundary conditions and the other handling the rest part of the domain. And in order to disentangle the two networks, a length factor function is introduced to eliminate the interference between them so as to make the method applicable to arbitrary geometries. We prove that the proposed method converges to the true solution of the boundary-value problem as the number of hidden units of the networks increases, and show by numerical experiments that PFNN can be applied to a series of linear and nonlinear second-order boundary value problems. Compared to existing approaches, the proposed method is able to produce more accurate solutions with fewer unknowns and less training costs.

The remainder of the paper is organized as follows. In Section 2, the basic framework of the PFNN method is presented. Following that we provide some theoretical analysis results on the accuracy of the PFNN method in Section 3. Some further comparisons between the present work and several recently proposed methods can be found in Section 4. We report numerical results on various second-order boundary value problems in Section 5. The paper is concluded in Section 6.

\section{The PFNN method}

Consider the following boundary-value problem:
\begin{equation}
  \label{BVP}
  \left \{
    \begin{array}{r l}
      -\nabla \cdot \Big( \rho(|\nabla u|) \nabla u \Big) + h(u) = 0,
      & \mbox{in}\ \Omega\subset \mathbb{R}^d, \\[2mm]
      u = \varphi,
      & \mbox{on}\ \Gamma_D, \\[1mm]
      \Big( \rho(|\nabla u|)\nabla u \Big)\cdot \bm{n} = \psi,
      & \mbox{on}\ \Gamma_N,
    \end{array}
  \right.
\end{equation}
where $\bm{n}$ is the outward unit normal,
$\Gamma_D\cup\Gamma_N = \partial\Omega$ and $\Gamma_D\cap\Gamma_N = \varnothing$.

It is well understood that under certain conditions (\ref{BVP}) can be seen as the Euler-Lagrange equation of the energy functional
\begin{equation}
  \label{Energy Functional}
  I[w]
  := \int_{\Omega} \left( P(w) + H(w) \right) d\bm{x}
   - \int_{\Gamma_N} \psi w d\bm{x},
\end{equation}
where
\begin{equation}
  \label{Energy Functional_P_H}
  P(w):= \int_{0}^{|\nabla w|} \rho(s)s ds
  \quad\mbox{and}\quad
  H(w):= \int_{0}^{w} h(s) ds.
\end{equation}

Introducing a hypothesis space $\mathcal{H}$ constructed with neural networks,
the approximate solution $u^{\ast}$ is obtained by solving the following minimization problem:
\begin{equation}
  \label{target}
  u^{\ast} = \arg\min\limits_{w\in \mathcal{H}} \Psi[w],
\end{equation}
where
\begin{equation}
  \label{Loss_Function}
  \Psi[w]
  := \displaystyle\frac{|\Omega|}{\#S(\Omega)}
     \sum\limits_{\bm{x}^i\in S(\Omega)} \left(
     P(w(\bm{x}^i)) +
     H(w(\bm{x}^i)) \right)
   - \displaystyle \frac{|\Gamma_N|}{\#S(\Gamma_N)}
     \sum\limits_{\bm{x}^i\in S(\Gamma_N)}
     \psi(\bm{x}^i) w(\bm{x}^i)
\end{equation}
is the loss function representing the discrete counterpart of $I[w]$,
$S(\Box)$ denotes a set of sampling points on $\Box$, and $\#S(\Box)$ is the size of $S(\Box)$.

In this work, the hypothesis space $\mathcal{H}$ is not an arbitrary space formed by neural networks.
Instead, it is designed to encapsulate the essential boundary conditions.
To construct the approximate solution $w_{\bm{\theta}}\in\mathcal{H}$,
we employ two neural networks, $g_{\bm{\theta_1}}$ and $f_{\bm{\theta_2}}$, instead of one, such that
\begin{equation}
  \label{TrialSolution}
  w_{\bm{\theta}}(\bm{x})  = g_{\bm{\theta_1}}(\bm{x}) + \ell(\bm{x}) f_{\bm{\theta_2}}(\bm{x}),
\end{equation}
where $\bm{\theta}=\{\bm{\theta_1},\bm{\theta_2}\}$ is the collection of the weights and biases of the two networks, and $\ell$ is a length factor function for measuring the distance to $\Gamma_D$, which satisfies
\begin{equation}
  \label{Ld_value}
  \left \{
    \begin{array}{l l}
      \ell(\bm{x})=0, & \bm{x}\in\Gamma_D, \\[1mm]
      \ell(\bm{x})>0, & \mbox{otherwise}.
    \end{array}
  \right.
\end{equation}

With the help of the length factor function, the neural networks $g_{\bm{\theta_1}}$ and $f_{\bm{\theta_2}}$ are utilized
to approximate the true solution $u$ on the essential boundary $\Gamma_D$ and
the rest of the domain, respectively.
The training of $g_{\bm{\theta_1}}$ and $f_{\bm{\theta_2}}$ are conducted separately,
where $g_{\bm{\theta_1}}$ is trained on the essential boundary to minimize the following energy functional
\begin{equation}
  \label{Energy_Functional_g}
    \Phi[g_{\bm{\theta_1}}] := \sum\limits_{\bm{x}^i\in S(\Gamma_D)}
    \left( \varphi(\bm{x}^i) - g_{\bm{\theta_1}}(\bm{x}^i) \right)^2,
\end{equation}
and  $f_{\bm{\theta_2}}$ is trained to approximate $u$ on the rest of the domain by minimizing the loss function (\ref{Loss_Function}) By this means, $f_{\bm{\theta_2}}$ would not produce any influence on $\Gamma_D$,
and the interference between the two networks are eliminated.

To construct the length factor function $\ell$,
the boundary of the domain is divided into $n$ segments: $\gamma_k$ ($k=1,2,\cdots,n$, $n\ge4$), with each segment either belonging to $\Gamma_D$ or to $\Gamma_N$. For any boundary segment  $\gamma_k\subset\Gamma_D$, select another segment $\gamma_{k_o}$ ($k_o\in\{1,2,\cdots,n\}$, $k_o\neq k$), which is not the neighbor of $\gamma_k$, as its companion. After that, for each $\gamma_k\subset\Gamma_D$, we construct a spline function $l_k$ that satisfies
\begin{equation}
  \label{Ld_basis}
  \left \{
    \begin{array}{l l}
      l_k(\bm{x}) = 0, & \bm{x}\in\gamma_k, \\[1mm]
      l_k(\bm{x}) = 1, &\bm{x}\in\gamma_{k_o},   \\[1mm]
      0<l_k(\bm{x})<1, & \mbox{otherwise}.
    \end{array}
  \right.
\end{equation}
Then $\ell$ is defined as:
\begin{equation}
  \label{Ld}
  \ell(\bm{x}) = \frac{\widetilde{\ell}(\bm{x})} {\max\limits_{\bm{x}\in\Omega} \widetilde{\ell}(\bm{x})},
  \quad\mbox{where}\quad
  \widetilde{\ell}(\bm{x}) = \prod\limits_{k\in\{k|\gamma_k\subset\Gamma_D\}} 1-(1-l_k(\bm{x}))^\mu.
\end{equation}
Here a hyper-parameter $\mu\geq 1$ is introduced to adjust the shape of $\ell$. A suggested value is $\mu = n_{\gamma_D}$, i.e., the number of the boundary segments on $\Gamma_D$, since in this way the average value of $\ell$ is kept in a proper range and would not decrease dramatically with the increase of $n_{\gamma_D}$.

For the purpose of flexibility, we employ the radial basis interpolation \cite{wright2003radial},
among other choices, to construct $l_k$.
Taking $m_k$ distinct points $\widehat{\bm{x}}^{k,1},\widehat{\bm{x}}^{k,2},\cdots,\widehat{\bm{x}}^{k,m_k}\in\gamma_k\cup \gamma_{k_o}$ as the interpolation nodes,  $l_k$ is defined as:
\begin{equation}
  \label{RBF_Intp}
  l_k(\bm{x}) = \sum_{i=1}^{m_k} a_{i} \phi(\bm{x};\widehat{\bm{x}}^{k,i})
         + \bm{b}\cdot\bm{x} + c,
\end{equation}
where 
\begin{equation}
  \label{RBF}
  \phi(\bm{x};\widehat{\bm{x}})
  = \left( e^2+\|\bm{x}-\widehat{\bm{x}}\|^2 \right)^{-1/2}
\end{equation}
is an inverse multiquadric radial basis function.
Here the parameter $e$ is used to adjust the shape of $\phi$,
To maintain scale invariance,
which is important for solving differential equations, we set {$e=1.25R/\sqrt{m_k}$},
where {$R$} is the {radius} of the minimal circle that encloses all the interpolation nodes.
Similar configurations of $e$ have been suggested in several previous studies
such as \cite{franke1982scattered, foley1987interpolation}.
The coefficients $a_{i},\bm{b}$ and $c$ are determined by solving a small linear system defined by the following constraints:
\begin{equation}
  \label{RBF_Intp_coef}
   \begin{array}{l}
      \sum\limits_{i=1}^{m_k} a_{i} \widehat{\bm{x}}^{k,i} = \bm{0},\quad
      \sum\limits_{i=1}^{m_k} a_{i} = 0,\quad
      l_k (\widehat{\bm{x}}^{k,i}) =  \left \{ \begin{array}{l}
      0, \,\,\, \widehat{\bm{x}}^{k,i}\in\gamma_k,\\[2mm]
      1, \,\,\, \widehat{\bm{x}}^{k,i}\in\gamma_{k_o},
      \end{array}\right.
      \,\,\, \forall i\in\{1,\cdots,m_k\}.
    \end{array}
\end{equation}
In particular, if $\gamma_k$ and $\gamma_{k_o}$ are two parallel hyper-planes, $l_k$ is reduced to a linear polynomial, which can be derived analytically without interpolation.

\begin{figure}[tb]
  \centering
  \subfigure[Boundary]
  {\includegraphics[width=0.24\textwidth]{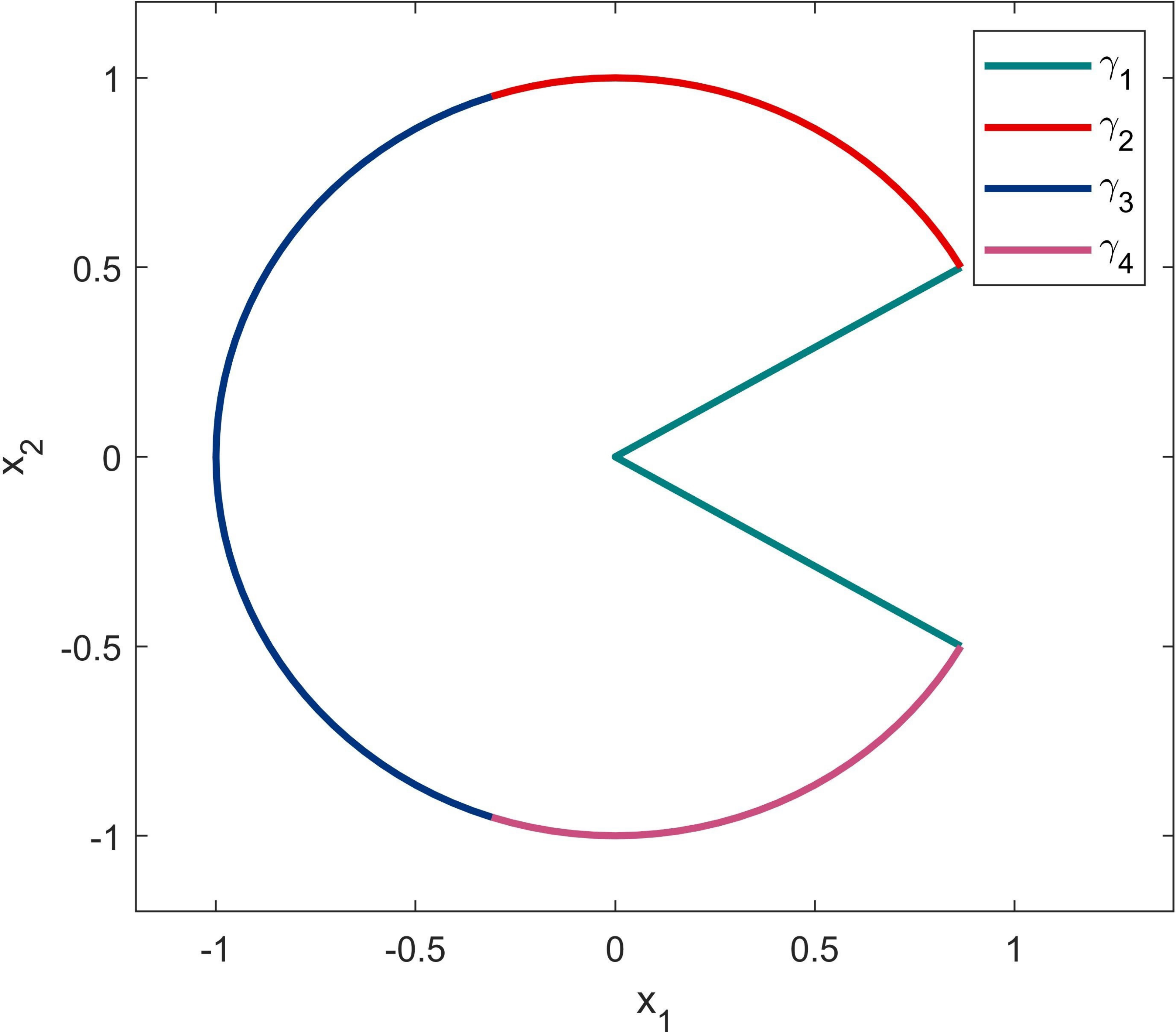}}
  \subfigure[$l_1$]
  {\includegraphics[width=0.24\textwidth]{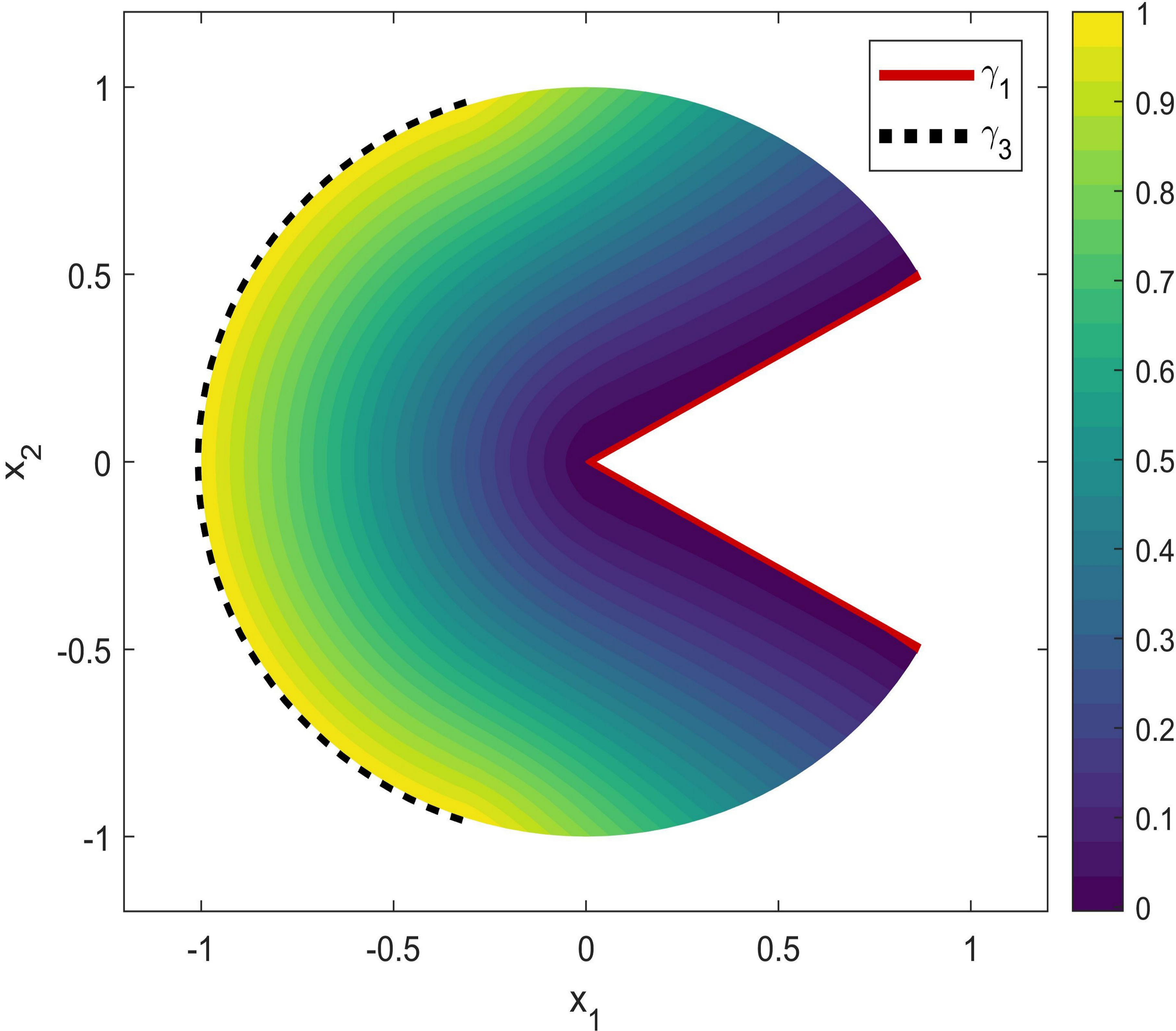}}
  \subfigure[$l_2$]
  {\includegraphics[width=0.24\textwidth]{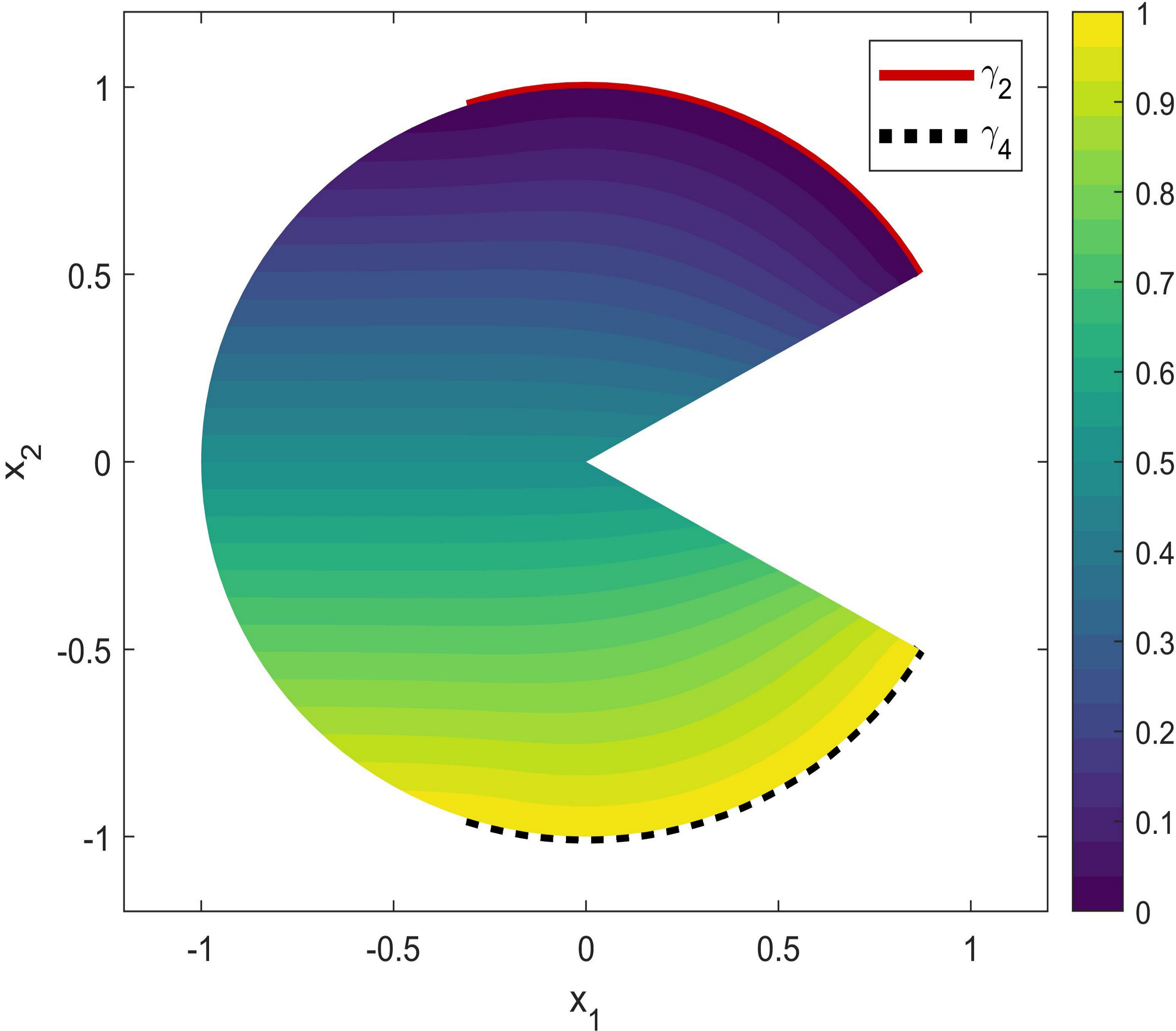}}
  \subfigure[$\ell$ ($\mu=2$)]
  {\includegraphics[width=0.24\textwidth]{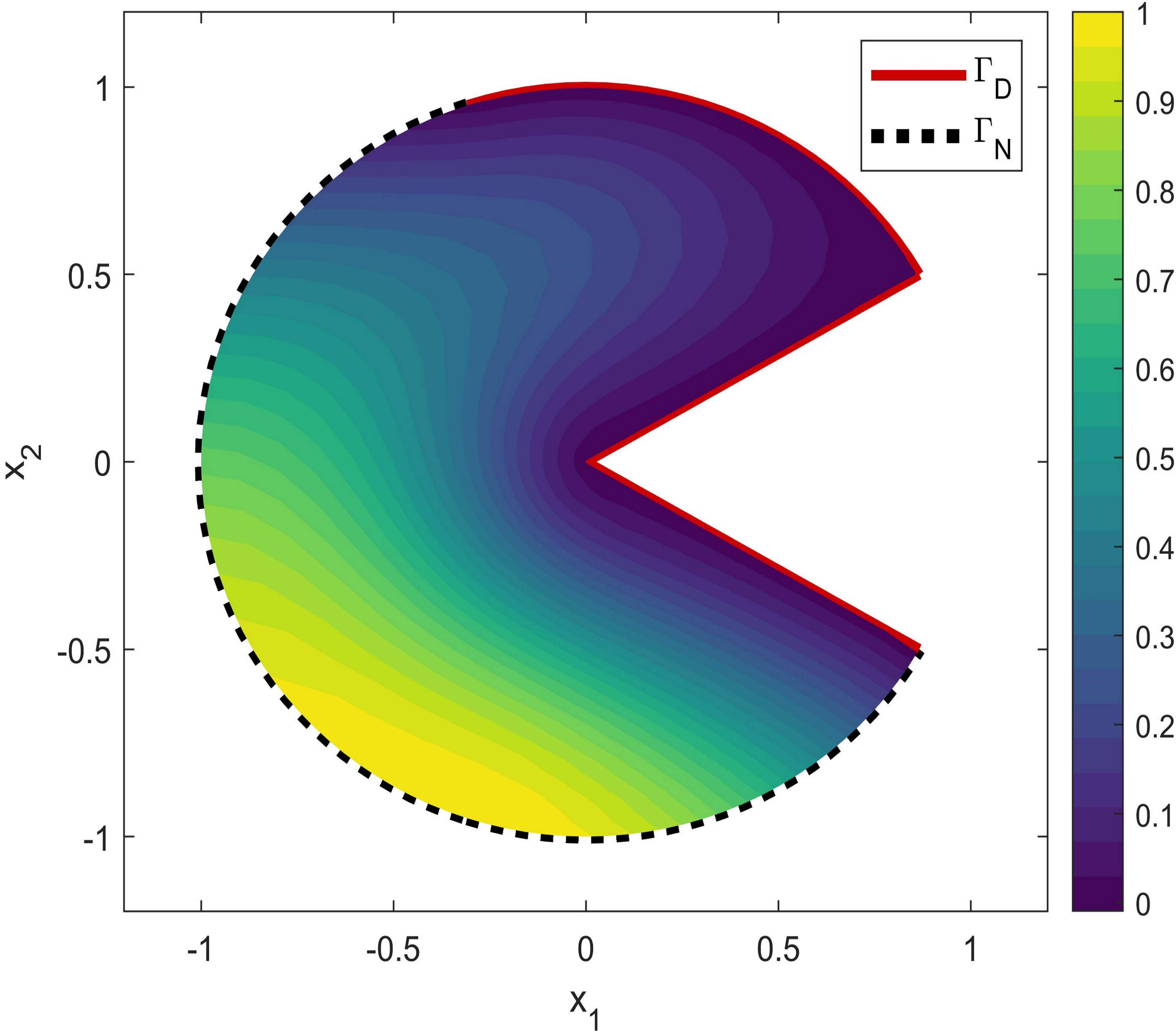}}
  \caption{An illustration on the construction of the length factor function $\ell$ for the case of a slotted disk.}
  \label{F_BVP_2D_Ld}
\end{figure}

To further illustrate how the length factor function $\ell$ is constructed,
consider an example that $\Omega$ is a slotted disk
as shown in Figure \ref{F_BVP_2D_Ld}.
Suppose the whole boundary $\partial\Omega=\Gamma_D\cup\Gamma_N$
is divided into four segments such that $\Gamma_D = \gamma_1\cup\gamma_2$
and {$\Gamma_N = \gamma_3\cup\gamma_4$},
where $\gamma_3$ and $\gamma_4$ are the opposite sides of
$\gamma_1$ and $\gamma_2$, respectively.
One can build $l_1$ and $l_2$ according to (\ref{RBF_Intp})-(\ref{RBF_Intp_coef})
and then construct $\ell$ with (\ref{Ld}).

We remark that such method is also feasible even when the boundary of the domain does not have an analytical form, so long as a set of sampling points $S(\partial\Omega)$ is available. Then by employing a certain clustering method one can divide $S(\partial\Omega)$ into $n$ subsets: $S(\gamma_k)$ ($k=1,2,\cdots,n$, $n\ge4$), with each subset responding to a boundary segment. Function $l_k$ can be constructed with the interpolation nodes taken from $S(\gamma_k\cup\gamma_{k_o})$ instead of $\gamma_k\cup\gamma_{k_o}$.

\section{Theoretical analysis of PFNN method}

In this section, we provide a theoretical proof to verify that under certain conditions the approximate solution obtained by the PFNN method is able to converge to the true solution of the boundary-value problem (\ref{BVP}) as the number of hidden units in the neural networks $g_{\bm{\theta_1}}$ and $f_{\bm{\theta_2}}$ increases.

We suppose that the true solution of (\ref{BVP}) belongs to the function space
\begin{equation}
  \label{Sobolev}
  \mathcal{W}^{1,p}_D(\Omega):=
  \{w\in \mathcal{W}^{1,p}(\Omega)|
    w=\varphi\ \mbox{on}\ \Gamma_D \},
\end{equation}
where
$\mathcal{W}^{1,p}(\Omega):= \{ w| \|w\|_{1,p,\Omega}<\infty \}$,
$\|w\|_{1,p,\Omega}:=
 \left( \displaystyle \int_{\Omega} \left(|w|^p + |\nabla w|^p\right) d\bm{x} \right)^{1/p}$.
Further, it is assumed that the function $\rho(\cdot)\cdot$ in (\ref{BVP}) satisfies the following conditions:
\begin{description}
  \item[(i)] $\rho(\cdot)\cdot$ is strictly increasing on $\mathbb{R}^{+}$ and $\rho(s)s|_{s=0}=0$.
  \item[(ii)] there exist constants $p>1$, $\lambda_1>0$ and $\lambda_2\geq0$, such that for all $s\in\mathbb{R}^{+}$, $\rho(s)s\leq \lambda_1 s^{p-1} - \lambda_2$.

  \item[(iii)] $\rho(\cdot)\cdot$ is H\"{o}lder continuous. If $1<p\leq2$, there exists a constant $C_1>0$, such that
      \begin{equation}
        \label{Continuous1}
        |\rho(s_1)s_1-\rho(s_2)s_2| \leq C_1|s_1-s_2|^{p-1},\ \forall s_1,s_2\in\mathbb{R}^{+};
      \end{equation}
      otherwise, there exist constants $C_1>0$, $K_1\geq0$ and $K_2>0$, such that
      \begin{equation}
        \label{Continuous2}
        |\rho(s_1)s_1-\rho(s_2)s_2| \leq C_1|s_1-s_2|(K_1+K_2(s_1+s_2)^{p-2}),\ \forall s_1,s_2\in\mathbb{R}^{+}.
      \end{equation}

  \item[(iv)] If $1<p\leq2$, there exist constants $C_2>0$, $K_1\geq0$ and $K_2>0$, such that
      \begin{equation}
        \label{Coercivity1}
        (\rho(s_1)s_1-\rho(s_2)s_2) (K_1+K_2(s_1+s_2)^{2-p})\geq C_2(s_1-s_2),\ \forall s_1\geq s_2\geq 0;
      \end{equation}
      otherwise, there exists a constant $C_2>0$, such that
      \begin{equation}
        \label{Coercivity2}
        \rho(s_1)s_1-\rho(s_2)s_2 \geq C_2(s_1^{p-1}-s_2^{p-1}),\ \forall s_1\geq s_2\geq 0.
      \end{equation}
\end{description}
Also, we assume that the function $h(\cdot)$ is monotonic increasing on $\mathbb{R}$ and satisfies
conditions similar to (\ref{Continuous1})-(\ref{Continuous2}).

We list two useful lemmas here for subsequent usage,
among which the second one is the famous universal approximation property of the neural networks.

\noindent\textbf{Lemma 1} (Chow, 1989, \cite{chow1989finite}) Suppose that the energy functional $I$ is in the form of (\ref{Energy Functional}),  $u=\arg\min\limits_{w\in\mathcal{W}^{1,p}_D(\Omega)} I[w]$, $u^{\ast}=\arg\min\limits_{w\in\mathcal{H}} I[w]$. Then there exists a constant $C>0$, such that
\begin{equation}
  \label{Lemma_Chow1}
  \left \{
    \begin{array}{l l}
      |u-u^{\ast}|_{1,p,\Omega}^2
      \leq C \|u-u^{\ast}\|_{1,p,\Omega}^{p-1}
      \inf\limits_{w\in\mathcal{H}}
      \|u-w\|_{1,p,\Omega}, & 1<p\leq 2, \\ [3mm]
      |u-u^{\ast}|_{1,p,\Omega}^p
      \leq C \|u-u^{\ast}\|_{1,p,\Omega}
      \inf\limits_{w\in\mathcal{H}}
      \|u-w\|_{1,p,\Omega}, & p>2.
    \end{array}
  \right.
\end{equation}
Especially, if $\mathcal{H}\subset\mathcal{W}^{1,p}_D(\Omega)$, then $u-u^{\ast} \in {W}^{1,p}_0(\Omega):= \{w| w\in\mathcal{W}^{1,p}(\Omega)\ \mbox{and}\ w=0,\ \mbox{on}\ \Gamma_D\}$. In this case, there exists a constant $C_3$ such that $\|u-u^{\ast}\|_{1,p,\Omega} \leq C_3|u-u^{\ast}|_{1,p,\Omega}$. Then we have
\begin{equation}
  \label{Lemma_Chow2}
  \|u-u^{\ast}\|_{1,p,\Omega}
  \leq \widetilde{C} \inf\limits_{w\in\mathcal{H}}
  \|u-w\|_{1,p,\Omega}^s,
\end{equation}
where $\widetilde{C}=C_3C$, $s=1/(3-p)$ if $1<p\leq2$ and $s=1/(p-1)$ if $p>2$.

\noindent\textbf{Lemma 2} (K. Hornik, 1991 \cite{hornik1991approximation}) Let
\begin{equation}
  \label{NN_SL}
  \mathcal{F}^{n}:=
  \{ f_{\bm{\theta}} |\
     f(\bm{x};\bm{\theta}) = \sum_{j=1}^{i} a_j\sigma(\bm{w}_j\cdot\bm{x}+b_j),\
     i\leq n \}
\end{equation}
be the space consisting of all neural networks with single hidden layer of no more than $n$ hidden units, where $\bm{\theta}=\{\bm{w}_j, a_j, b_j, j=1,2,\cdots,i\}$, $\bm{w}_j$ is the input weight, $a_j$ is the output weight and $b_j$ is the bias. If $\sigma\in\mathcal{W}^{m,p}(\Omega)$ is nonconstant and all its derivatives up to order $m$ are bounded, then for all $f\in\mathcal{W}^{m,p}(\Omega)$ and $\varepsilon>0$, there exists an integer $n>0$ and a function $f_{\bm{\widehat{\theta}}} \in\mathcal{F}^n$, such that
\begin{equation}
  \label{Universal_Approximation}
  \|f-f_{\bm{\widehat{\theta}}}\|_{m,p,\Omega} < \varepsilon.
\end{equation}

The main convergence theorem of the PFNN method is given below.

\noindent\textbf{Theorem 1} Suppose that the energy functional $I$ is in the form of (\ref{Energy Functional}), $u = \arg\min\limits_{w\in\mathcal{W}^{1,p}_D(\Omega)} I[w]$,
$g^{n_1} = \arg\min\limits_{g_{\bm{\theta_1}}\in\mathcal{F}^{n_1}}
 \displaystyle\int_{\Gamma_D}
 (\varphi-g_{\bm{\theta_1}})^2 dx$,
$u^{n_1,n_2} =
\arg\min\limits_{w_{\bm{\theta}}\in\mathcal{H}^{n_1,n_2}}
I[w_{\bm{\theta}}]$,
where
$\mathcal{H}^{n_1,n_2}:= \{w_{\bm{\theta}}|
w_{\bm{\theta}} = g^{n_1} + \ell f_{\bm{\theta_2}},\
f_{\bm{\theta_2}}\in\mathcal{F}^{n_2}\}$,
then
\begin{equation}
  \label{Method_Error}
  \|u-u^{n_1,n_2}\|_{1,p,\Omega} \rightarrow 0
  \quad \mbox{as}\quad
  n_1\rightarrow\infty
  \ \mbox{and}\
  n_2\rightarrow\infty.
\end{equation}

\noindent\emph{\textbf{Proof}} For the sake of brevity, we drop the subscripts in $w_{\bm{\theta}}$, $g_{\bm{\theta_1}}$ and $f_{\bm{\theta_2}}$ and utilize $w$, $g$ and $f$ to represent them, respectively. We only need to consider the case of $1<p\leq2$ since the proof for the case of $p>2$ is analogous. Combining the Friedrichs' inequality with {Lemma 1}, we have
\begin{equation}
  \label{Approxiamte_Error}
  \begin{array}{l l}
    \|u-u^{n_1,n_2}\|_{1,p,\Omega}
    & \leq C_4 \left( |u-u^{n_1,n_2}|_{1,p,\Omega} +
                      \|u-u^{n_1,n_2}\|_{0,{2},\Gamma_D}\right) \\ [1mm]
    & \leq C_5 \left( \|u-u^{n_1,n_2}\|_{1,p,\Omega}^{(p-1)/2}
                      \inf\limits_{w\in\mathcal{H}^{n_1,n_2}}
                      \|u-w\|_{1,p,\Omega}^{1/2} +
                      \|u-u^{n_1,n_2}\|_{0,2,\Gamma_D} \right),
  \end{array}
\end{equation}
where $C_4$, $C_5$ are constants.

First, we prove that $\|u-u^{n_1,n_2}\|_{0,2,\Gamma_D}\rightarrow 0$ as $n_1\rightarrow\infty$, for arbitrary $n_2\in\mathbb{N}^{+}$. According to Lemma 2, for all $\varepsilon>0$, there exists an integer $n_1^{\ast}>0$ and a function $\widehat{g} \in\mathcal{F}^{n_1^{\ast}}$ satisfying $\|\varphi-\widehat{g}\|_{0,2,\Gamma_D}<\varepsilon$.
Due to
$g^{n_1} = \arg\min\limits_{g\in\mathcal{F}^{n_1}}
 \displaystyle\int_{\Gamma_D} (\varphi-g)^2 dx
 = \arg\min\limits_{g\in\mathcal{F}^{n_1}}
 \| \varphi-g \|_{0,2,\Gamma_D}$,
the following relationship hold:
\begin{equation}
  \label{Boundary_Errors}
  \|u-u^{n_1,n_2}\|_{0,2,\Gamma_D}
  = \|\varphi-g^{n_1}\|_{0,2,\Gamma_D}
  \leq \|\varphi-\widehat{g}\|_{0,2,\Gamma_D}
  < \varepsilon,\quad
  \forall\ n_1\geq n_1^{\ast}.
\end{equation}

We then prove that $\inf\limits_{w\in\mathcal{H}^{n_1,n_2}} \|u-w\|_{1,p,\Omega}^{1/2}\rightarrow 0$ as $n_2\rightarrow\infty$, for arbitrary $n_1\in\mathbb{N}^{+}$.
Since $u$, $g^{n_1}$, $\ell\in\mathcal{W}^{1,p}(\Omega)$ and $\ell>0$, in $\Omega$, we have $\displaystyle\frac{u-g^{n_1}}{\ell}\in\mathcal{W}^{1,p}(\Omega)$. According to Lemma 2, for all $\varepsilon>0$ and $n_1>0$, there exists an integer $n_2^{\ast}>0$, which is dependent on $g^{n_1}$ and therefore relies on $n_1$, and a function $\widehat{f} \in\mathcal{F}^{n_2^{\ast}}$, such that
\begin{equation}
  \| \displaystyle\frac{u-g^{n_1}}{\ell} -
     \widehat{f} \|_{1,p,\Omega}
  < \frac{\varepsilon^2}{\|\ell\|_{1,p,\Omega}}.
\end{equation}
Correspondingly, the function $\widehat{w} = g^{n_1} + \ell \widehat{f} \in\mathcal{H}^{n_1,n_2^{\ast}}$ satisfies
\begin{equation}
  \|u-\widehat{w}\|_{1,p,\Omega}
  = \| \ell(\frac{u-g^{n_1}}{\ell} -
       \widehat{f}) \|_{1,p,\Omega}
  \leq \|\ell\|_{1,p,\Omega}
       \| \frac{u-g^{n_1}}{\ell} -
          \widehat{f} \|_{1,p,\Omega}
  < \varepsilon^2.
\end{equation}
It then follows that
\begin{equation}
  \label{Space_Error}
  \inf\limits_{w\in\mathcal{H}^{n_1,n_2}}
  \|u-w\|_{1,p,\Omega}^{1/2}
  \leq \|u-\widehat{w}\|_{1,p,\Omega}^{1/2}
  < \varepsilon,\quad
  \forall n_2\geq n_2^{\ast}(n_1).
\end{equation}

Finally, we make use of reduction to absurdity to prove that $\|u-u^{n_1,n_2}\|_{1,p,\Omega}$ can be arbitrarily small so long as $n_1$ and $n_2$ are large enough. Suppose that for all $n_1>0$ and $n_2>0$, there exists a constant $\delta>0$, such that $\|u-u^{n_1,n_2}\|_{1,p,\Omega} \geq \delta$. According to (\ref{Boundary_Errors}) and (\ref{Space_Error}), there exist integers $\widehat{n}_1>0$ and $\widehat{n}_2>0$ ($\widehat{n}_2$ is dependent on $\widehat{n}_1$), such that $\|u-u^{\widehat{n}_1,\widehat{n}_2}\|_{0,2,\Gamma_D}
 < \displaystyle\frac{\delta}{2C_5}$ and
$\inf\limits_{w \in\mathcal{H}^{\widehat{n}_1,\widehat{n}_2}}
 \|u-w\|_{1,p,\Omega}^{1/2}
 < \displaystyle\frac{\delta^{(3-p)/2}}{2C_5}$.
Combining with (\ref{Approxiamte_Error}), we have
\begin{equation}
  \label{Contradiction_1}
  \begin{array}{l}
    \begin{array}{l l}
      \|u-u^{\widehat{n}_1,\widehat{n}_2}\|_{1,p,\Omega}^{(3-p)/2}
      & \leq C_5 \left(
        \inf\limits_{w \in\mathcal{H}^{\widehat{n}_1,\widehat{n}_2}}
        \|u-w\|_{1,p,\Omega}^{1/2} +
        \|u-u^{\widehat{n}_1,\widehat{n}_2}\|_{0,2,\Gamma_D}
        \|u-u^{\widehat{n}_1,\widehat{n}_2}\|_{1,p,\Omega}^{(1-p)/2}
        \right) \\
      & \leq C_5 \left(
        \inf\limits_{w \in\mathcal{H}^{\widehat{n}_1,\widehat{n}_2}}
        \|u-w\|_{1,p,\Omega}^{1/2} +
        \|u-u^{\widehat{n}_1,\widehat{n}_2}\|_{0,2,\Gamma_D}
        \delta^{(1-p)/2}
        \right) \\
      & < C_5 \left(
        \displaystyle\frac{\delta^{(3-p)/2}}{2C_5} +
        \displaystyle\frac{\delta}{2C_5}
        \delta^{(1-p)/2}
        \right)
        = \delta^{(3-p)/2},
    \end{array}
  \end{array}
\end{equation}
which is contradictory to $\|u-u^{\widehat{n}_1,\widehat{n}_2}\|_{1,p,\Omega} \geq \delta$. Therefore we can conclude that $\|u-u^{n_1,n_2}\|_{1,p,\Omega} \rightarrow 0$ as $n_1\rightarrow\infty$ and $n_2\rightarrow\infty$.
$\hfill\square$ 

We remark there that although only neural networks with single hidden layers are considered here,  similar conclusions can be made for the case of multi-layer neural networks.
The details are omitted for brevity.
It is also worth noting that the assumptions in the analysis are only sufficient conditions
for the converge of the PFNN method. Later numerical experiments
with cases that the assumptions are not fully satisfied will show that the proposed PFNN method still works well.

\section{Comparison with other methods}

Most of the existing neural network methodologies employ a penalty method to deal with the essential boundary conditions. For the boundary-value problem (\ref{BVP}), a straightforward way is to minimize the following energy functional in least-squares form:
\begin{equation}
  \label{Loss_Least_Square}
  \begin{array}{l l}
    {I}[w]
    & := \displaystyle\int_{\Omega} \Big(
         - \nabla \cdot ( \rho(|\nabla w|) \nabla w )
         + h(w) \Big)^2 d\bm{x}
       + \beta_1 \displaystyle\int_{\Gamma_D}
         (\varphi-w)^2 d\bm{x} \\
    &  +\ \beta_2 \displaystyle\int_{\Gamma_N} \Big(
         ( \rho(|\nabla w|)\nabla w )\cdot \bm{n}
         - \psi\Big)^2 d\bm{x},
  \end{array}
\end{equation}
where $\beta_1$ and $\beta_2$ are penalty coefficients, which can also be seen as Lagrangian multipliers
that help transform a constrained optimization problem into an unconstrained one.

Such least-squares approach has not only been adopted in many early studies
\cite{van1995neural, lagaris1998artificial, mcfall2009artificial, jianyu2002numerical, mai2001numerical, li2010integration},
but also employed by a number of recent works based on deep neural networks.
An example is the physics-informed neural networks \cite{raissi2019physics} in which
deep neural networks are applied to solve forward and inverse problems involving nonlinear partial differential equations
arising from thermodynamic, fluid mechanics and quantum mechanics.
Following it, in the work of the hidden fluid mechanics \cite{raissi2020hidden}
impressive results are obtained by extracting velocity and pressure information from the data of flow visualizations.
Another example of using least-squares energy functionals is the work of the Deep Galerkin method \cite{sirignano2018dgm},
in which efforts are made to merge deep learning techniques with the Galerkin method for
efficiently solving several high-dimensional differential equations.
Overall, a major difficulty by using the least-squares approach is that approximations have to be made
to the high-order derivatives of the true solution in some way, which could eventually
lead to high  training cost and low accuracy.

To avoid approximating high-order derivatives, several methods have been proposed to transform
the original problem into the corresponding weak forms.
Examples include the Deep Ritz method \cite{weinan2018deep}
which employs an energy functional of the Ritz type:
\begin{equation}
  \label{Loss_DRM}
  {I}[w]
  := \int_{\Omega} \left( P(w) + H(w) \right) d\bm{x}
   - \int_{\Gamma_N} \psi w d\bm{x}
   + \beta \int_{\Gamma_D} (\varphi-w)^2 d\bm{x},
\end{equation}
and the Deep Nitsche method \cite{liao2019deep} which is based on an energy functional in the sense of Nitsche:
\begin{equation}
  \label{Loss_DNM}
  \begin{array}{l l}
    {I}[w]:=
    & \displaystyle\int_{\Omega} \left( P(w) + H(w) \right) d\bm{x}
      - \int_{\Gamma_N} \psi w d\bm{x}
      - \displaystyle\int_{\Gamma_D}
        w \Big(\rho(|\nabla w|) \nabla w \cdot \bm{n}\Big) d\bm{x} \\
    & - \displaystyle\int_{\Gamma_D}
         \varphi \left( \beta w - \Big(\rho(|\nabla w|)\nabla w \cdot \bm{n}\Big)
             \right) d\bm{x}
      + \frac{\beta}{2}
        \displaystyle\int_{\Gamma_D} w^2 d\bm{x}.
  \end{array}
\end{equation}
The transformations to weak forms done in the two methods can effectively reduce the smoothness requirement
of the approximate solution. However, the penalty terms due to the essential boundary condition still persist,
which would lead to extra training cost. Moreover, there is little clue on how to set the penalty factor
$\beta$; improper values will cause negative influence on the accuracy of the approximate solution
or even lead to training failure.

There are also several attempts \cite{lagaris1998artificial, lagaris2000neural, mcfall2009artificial, mcfall2013automated} made to eliminate the penalty terms by explicitly manufacturing an approximate solution in the following form:
\begin{equation}
  \label{TrialSolution2}
  w_{\bm{\theta}}(\bm{x}) =
  \mathcal{G}(\bm{x}) + \mathcal{L}(\bm{x}) f_{\bm{\theta}}(\bm{x}),
\end{equation}
where $\mathcal{G}$ satisfies the essential boundary condition and
$\mathcal{L}$ serves as a length factor function.
It is worth noting that  the scopes of applications of these methods are rather narrow,
despite that the approximate solutions do share  some similarities with the PFNN method.
In particular, these methods usually construct function $\mathcal{G}$
either in analytic forms \cite{lagaris1998artificial, lagaris2000neural},
or through spline interpolations \cite{mcfall2009artificial, mcfall2013automated},
and are therefore only suitable to simple geometries in low dimensions.
To establish the length factor function $\mathcal{L}$,
these methods usually rely on mapping the original domain to a hyper-sphere,
which is again not flexible and efficient for problems with complex geometries and high dimensions.

The proposed PFNN method can effectively combine the advantages of the aforementioned state-of-the-art
while overcoming their drawbacks. It reduces the smoothness requirements as well as removes
the penalty term in the loss function, effectively converting the original hard problem to a relatively easy one.
By introducing two neural networks, instead of only one, the approximations made to the true solutions
are separated to the essential boundaries and the rest of the domain, respectively.
The original training task is divided into two simpler ones, which can substantially reduce
the training cost as well as enhance the accuracy of the approximation.
To eliminate the interference between the two neural networks, a systematic approach is
further proposed to construct the length factor function in a most flexible manner.
As we will show later in the numerical experiments,
the PFNN method is applicable to a wide range of problems on complex geometries in arbitrary dimensions
and is able to achieve higher accuracy with lower training cost as compared with other approaches.

\section{Numerical experiments}

\begin{table}[!htb]
  \caption{List of test cases for numerical experiments.}
  \label{Test_case_design}
  \centering
  \small
  \renewcommand{\arraystretch}{1.2}
  \begin{tabu}{c|c|c|c|c|c}
    \tabucline[1pt]{-}
    Case & Problem               & Dimension & Domain   & $\rho(|\nabla u|)$ & $h(u)$ \\
    \hline
    1    & anisotropic diffusion & 2   & square         & $\bm{A}$ & c \\
    2    & minimal surface       & 2   & Koch Snowflake & $1/\sqrt{1+|\nabla u|^2} $ & 0 \\
    3    & $p$-Liouville-Bratu   & 3   & Stanford Bunny & $|\nabla u|^{p-2}$ & $-\lambda\exp(u)+c$ \\
    4    & Poisson-like    & 100 & hypercube      & $|\nabla u|^{p-2}$ & $u+c$ \\
    \tabucline[1pt]{-}
  \end{tabu}
\end{table}

\begin{figure}[!htb]
  \centering
  \subfigure[Koch Snowflake ($L=5$)]
  {\includegraphics[height=6.0cm]{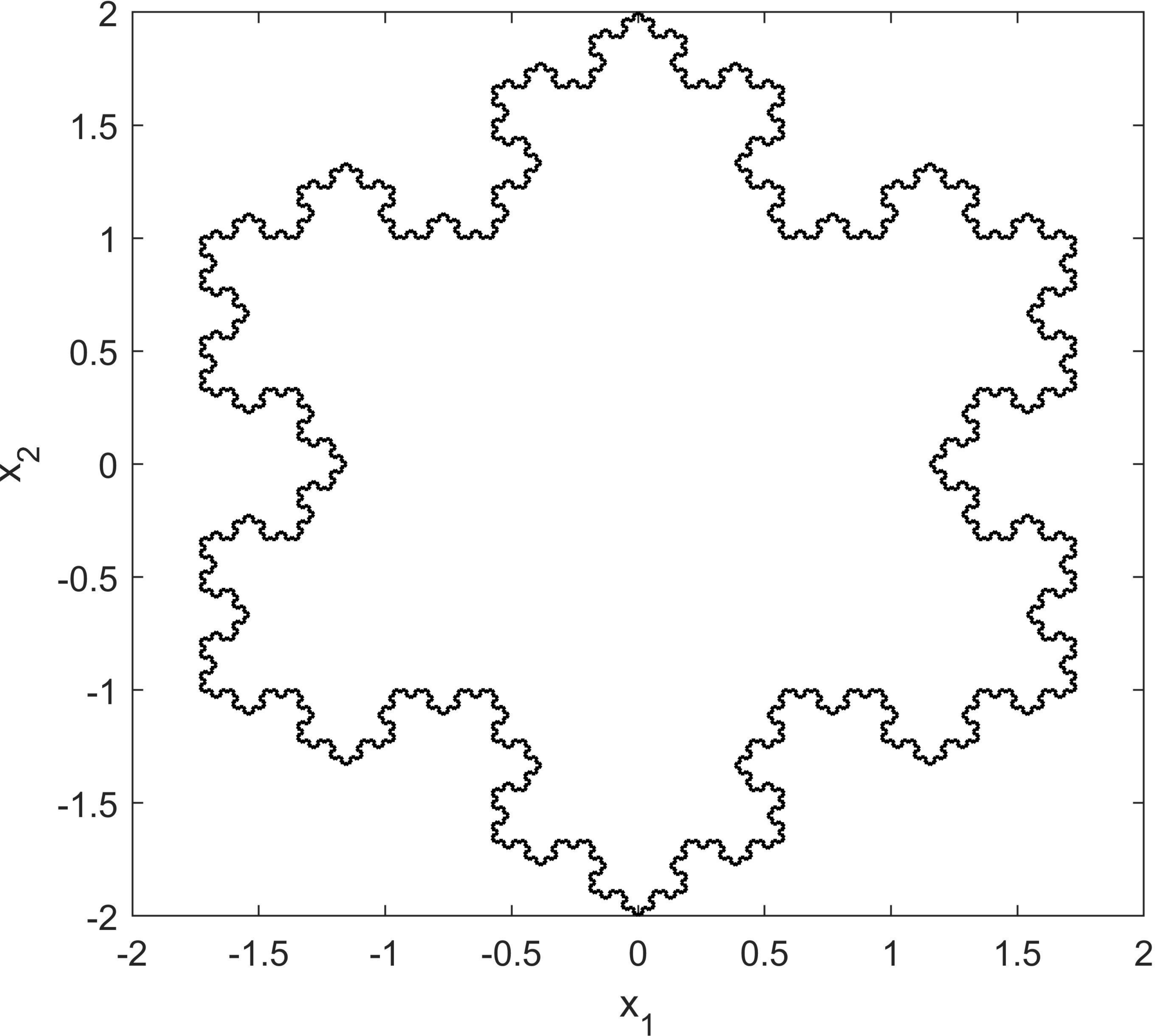}}
  \subfigure[Stanford Bunny]
  {\includegraphics[height=5.8cm]{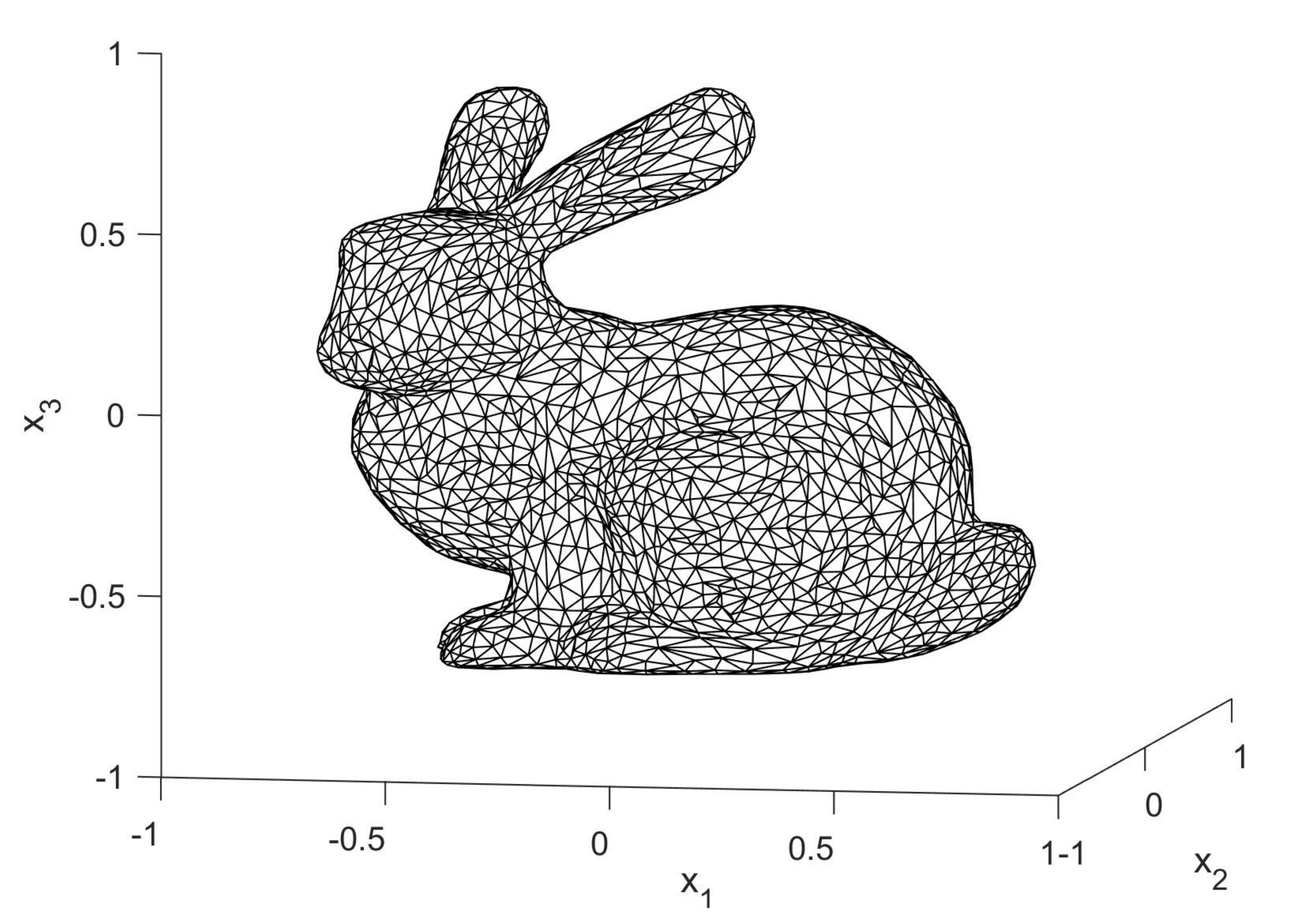}}
  \caption{The computation domains with complex geometries in test cases 2 (a) and 3 (b).}
  \label{fig:Test_case_domain}
\end{figure}

A series of numerical experiments are conducted to examine the numerical behaviors of the proposed PFNN method
as well as several previous state-of-the-art approaches.
We design four test cases covering different types of problems that have the same form of (\ref{BVP})
but vary in $\rho(|\nabla u|)$ and $h(u)$, as shown in Table \ref{Test_case_design}.
In particular, the computational domains of test cases 2 and 3 are further illustrated
in Figure \ref{fig:Test_case_domain}, both exhibiting very complex geometries.
We employ the ResNet model \cite{he2016deep} with sinusoid activation functions
to build the neural networks in both PFNN and other neural network based methods.
The Adam optimizer \cite{kingma2014adam} is utilized for training,
with the initial learning rate set to 0.01.
Unless mentioned otherwise, the maximum number of iteration is set to 5,000 epochs.
The relative $\ell_2$-norm
is used to estimate the error of the approximate solution. 
In all tests except those with the traditional finite element method,
we perform ten independent tests and collect the results for further analysis.

\subsection{Anisotropic diffusion equation on a square}
The first test case is an anisotropic diffusion equation:
\begin{equation}
  \label{PDE_Poisson}
  -\nabla\cdot\bm{A}\nabla u + {c} = 0\quad \mbox{in}\,\,\Omega=[-1,1]^2,
\end{equation}
where \begin{equation}
  \bm{A} = \left[
  \begin{array}{l l}
    (x_1+x_2)^2+1  & {-x_1^2+x_2^2} \\[1mm]
    {-x_1^2+x_2^2} & (x_1-x_2)^2+1  \\
  \end{array}
\right].
\end{equation}
The corresponding energy functional to minimize is:
\begin{equation}
  \label{EF_Poisson}
  I[w] = \int_{\Omega}
  ({\frac{1}{2}}\bm{A}\nabla w \cdot \nabla w + {c} w ) d\bm{x}
   - \int_{\Gamma_N} \psi w d\bm{x}.
\end{equation}
In the experiment, we set the essential boundary to $\Gamma_D = \{-1\}\times[-1,1]$
and set the exact solution to various forms including
$u^{(1)} = \ln(10(x_1+x_2)^2 + (x_1-x_2)^2 + 0.5)$,
$u^{(2)} = (x_1^3-x_1)\cosh(2x_2)$,
$u^{(3)} = (x_1^2-x_2^2) / (x_1^2+x_2^2+0.1)$, and
$u^{(4)} = (\sign x_1) \left((x_1+\sign x_1)^4 - 1\right) \exp(-x_2^2)$,
where $u^{(4)}$, discontinuous as such, is a weak solution of the partial differential equation.

For comparison purpose, we examine the performance of several approaches including the linear finite element method,
the least-squares neural network method, the Deep Ritz method, the Deep Nitsche method
and our proposed PFNN method.
A uniform mesh of $30\times30$ cells is used by the bilinear finite element method,
leading to $31\times30=930$ unknowns.
For the penalty-based neural network methods, the network adopted is comprised of 4 ResNet
blocks, each of which contains 2 fully connected layers with 10 units and a residual connection,
resulted in totally $811$ undecided parameters.
And for the PFNN method,  the network $g_{\bm{\theta_1}}$ and $f_{\bm{\theta_2}}$ consist of 1
and 3 ResNet blocks, respectively, corresponding to a total of  $151+591=742$ parameters.
The penalty coefficients in the penalty-based approaches are all
set to three typical values: $100$, $300$ and $500$.

\begin{figure}[!htb]
  \centering
  \subfigure[$u^{(1)} = \ln(10(x_1+x_2)^2 + (x_1-x_2)^2 + 0.5)$]
  {\includegraphics[width=0.49\textwidth]{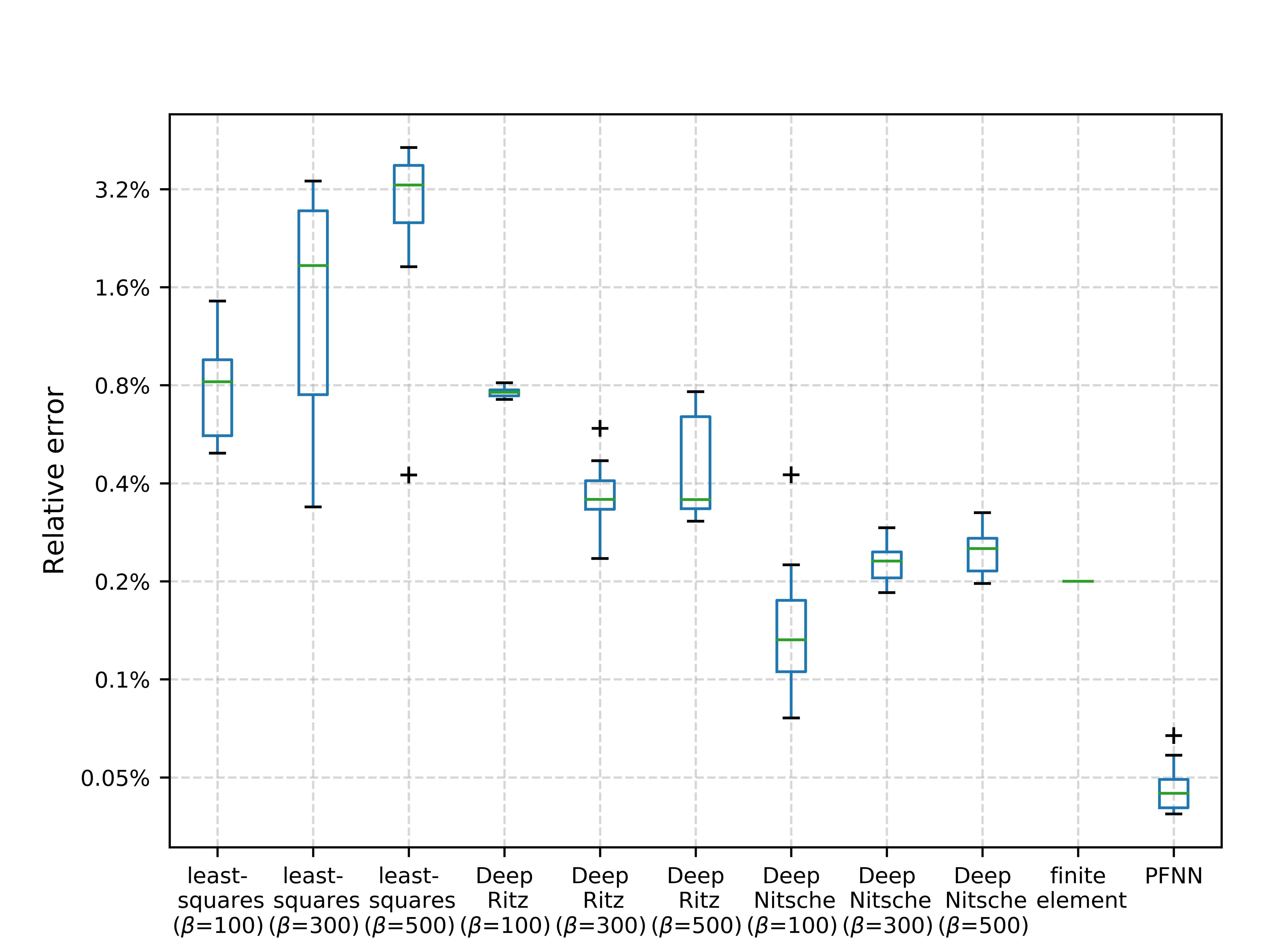}}
  \subfigure[$u^{(2)} = (x_1^3-x_1)\cosh(2x_2)$]
  {\includegraphics[width=0.49\textwidth]{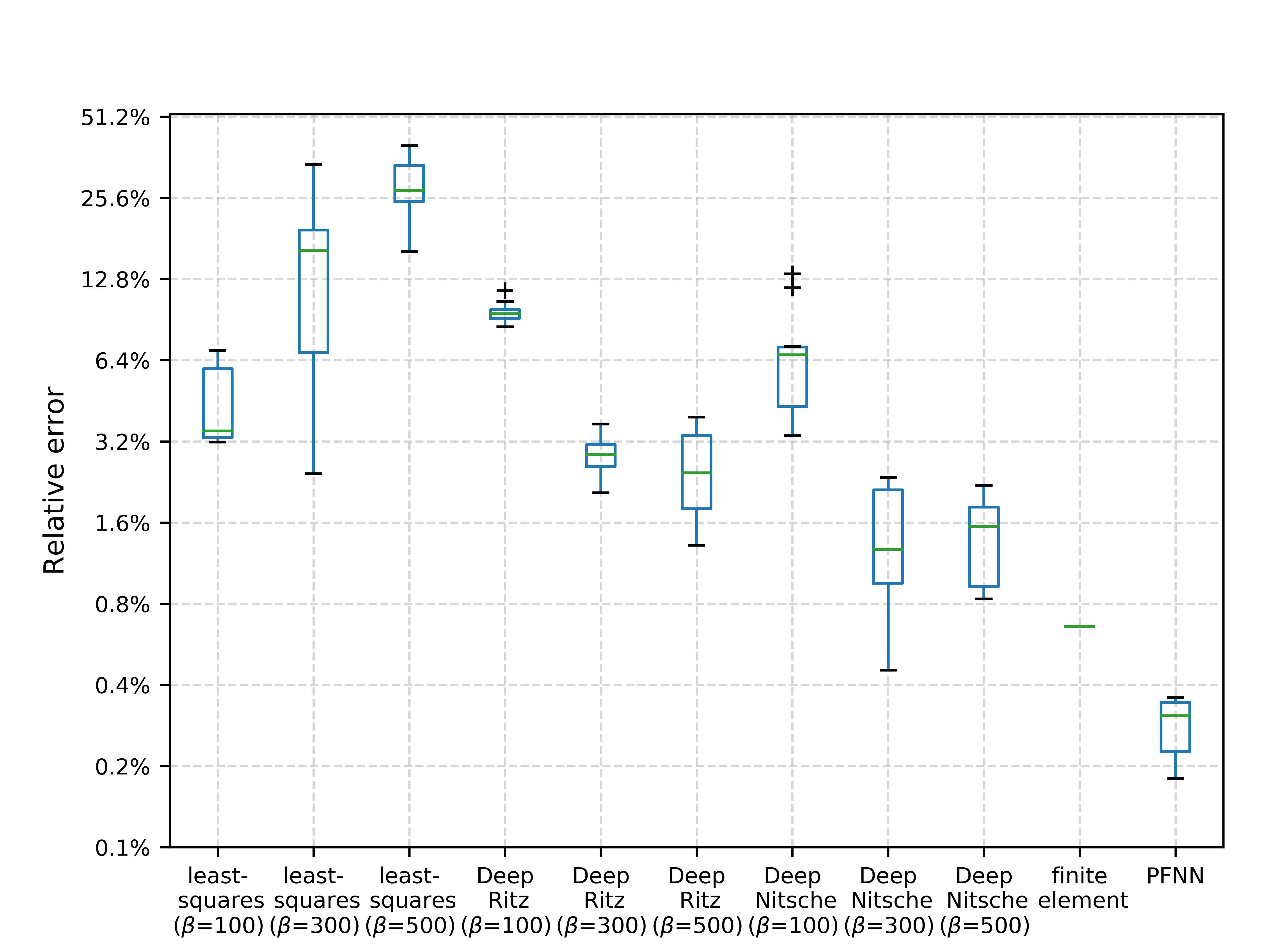}} \\
  \subfigure[$u^{(3)} = (x_1^2-x_2^2)/(x_1^2+x_2^2+0.1)$]
  {\includegraphics[width=0.49\textwidth]{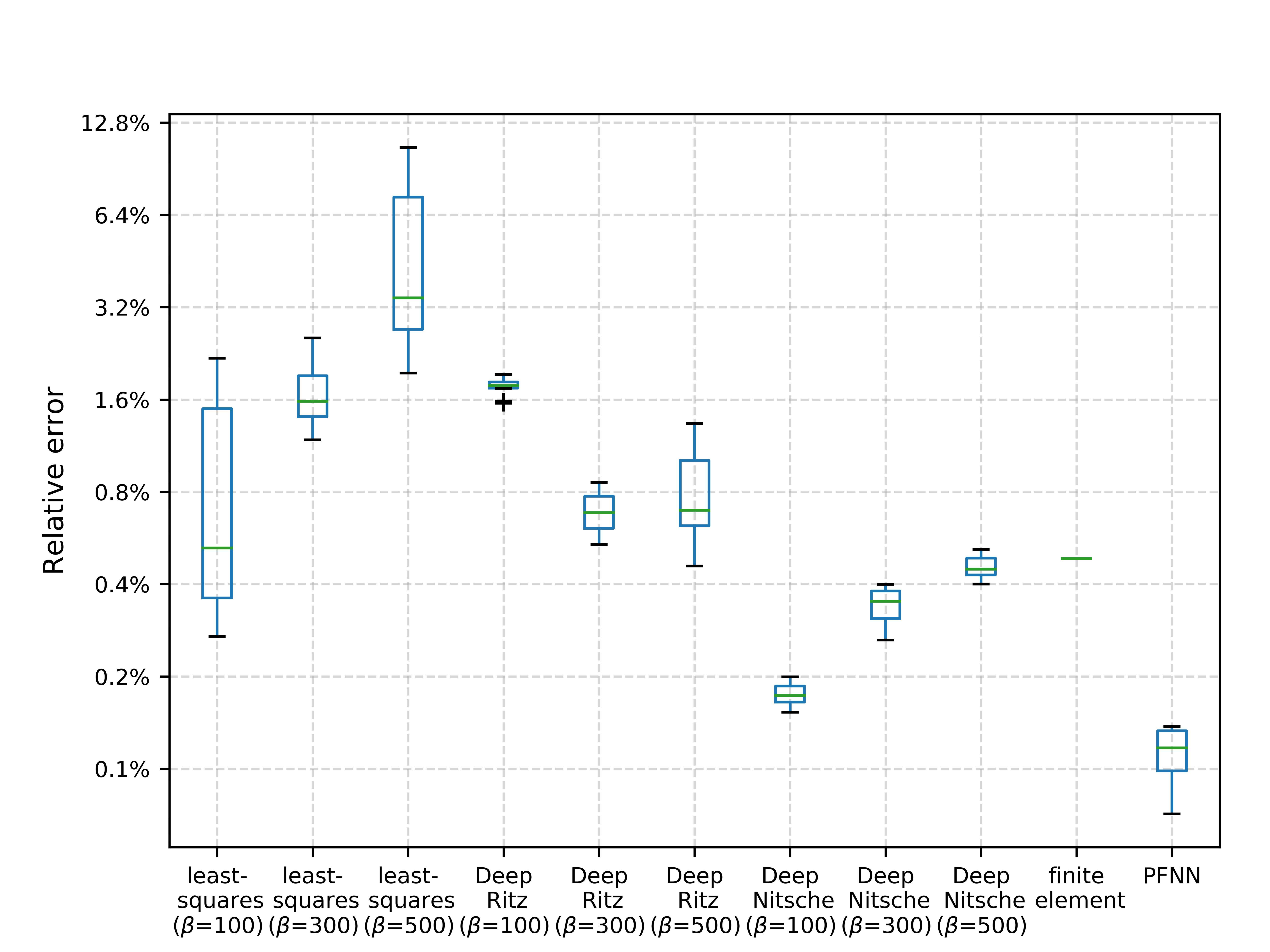}}
  \subfigure[$u^{(4)} = \mbox{sign}(x_1) \Big((x_1+\mbox{sign}(x_1))^4 - 1\Big) \exp(-x_2^2)$]
  {\includegraphics[width=0.49\textwidth]{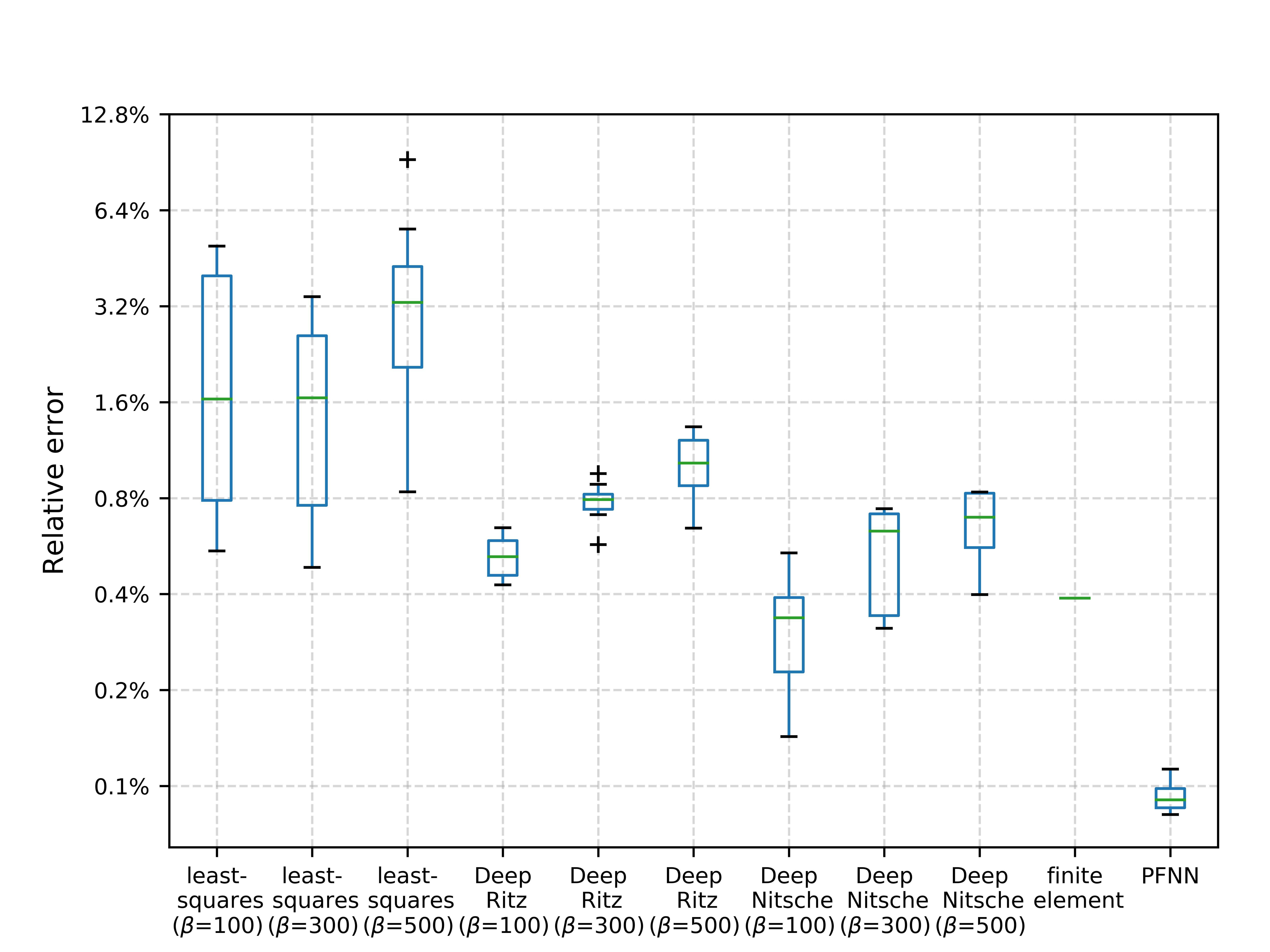}}
  \caption{Results of solving the anisotropic diffusion equation on a square with various methods.}
  \label{F_BVP_2D_Poisson_Rectangle}
\end{figure}

The experiment results are illustrated in Figure \ref{F_BVP_2D_Poisson_Rectangle}.
To get detailed information of the achievable accuracy of every approach,
we draw in the figure  the box-plots that contain {the five-number summary}
including the smallest observation, lower quartile, median, upper quartile, and largest observation for each set of tests.
From the figure, it can be observed that the performance of the classical least-squares neural network method
is usually the worst, due in large part to the approximations made to high-oder derivatives.
By introducing the corresponding weak forms, Deep Ritz and Deep Nitsche methods
can deliver better results than the least-squares method does,
but the performance still has a strong dependency on the specific values of the penalty coefficients
and the performance is in many cases not competitive to the traditional finite element method.
The advantages of the PFNN method are quite clear:
it is able to outperform all above approaches in all tested problems
in terms of the sustained accuracy and is much more robust than the penalty-based methods.

\subsection{Minimal surface equation on a Koch snowflake}
Consider a minimal surface equation \cite{giusti1984minimal}:
\begin{equation}
  \label{PDE_MinSurface}
  -\nabla \cdot \left(\frac{\nabla u}{\sqrt{1+|\nabla u|^2}} \right) = 0,
\end{equation}
defined on a well-known fractal pattern -- the Koch snowflake domain \cite{koch1904courbe},
as shown in Figure \ref{fig:Test_case_domain} (a).
The equation (\ref{PDE_MinSurface}) can be seen as the Euler-Lagrange equation of the energy functional:
\begin{equation}
  \label{EF_MinSurface}
  I[w]:= \int_{\Omega} \sqrt{1+|\nabla w|^2} d\bm{x},
\end{equation}
which represents the area of the surface of $w$ on $\Omega$.

The minimal surface problem is to find a surface with minimal area under the given boundary conditions. It can be verified that the catenoid surface
\begin{equation}
  \label{Catenoid}
  u = \lambda\cosh(\displaystyle\frac{x_1}{\lambda}) \sqrt{1-z^2},\quad z=\left(\displaystyle\frac{x_2}{\lambda\cosh(x_1/\lambda)}\right)
\end{equation}
satisfies the equation (\ref{PDE_MinSurface}) for all $\lambda$ satisfying $\left|z\right| \leq 1$ in $\Omega$.
In particular, for the Koch snowflake, the condition becomes $\lambda\geq 2$.
The difficulty of the problem is increased as $\lambda\rightarrow 2$
(especially, $|\nabla u|$ is unbounded on $\Omega$ when $\lambda=2$).
To examine how well various methods can deal with such kind of problem,
we set $\lambda$ to a relatively challenging value: $\lambda=2.1$.

\begin{table}[!htb]
  \caption{Errors of various methods for solving the minimal surface equation on a Koch snowflake.}
  \label{T_BVP_2D_MinSurface_Snow}
  \centering
  \small
  \renewcommand{\arraystretch}{1.0}
  \begin{tabu}{c|l l|l l|c}
    \tabucline[1pt]{-}
          Method   & \multicolumn{2}{c|}{Deep Ritz} & \multicolumn{2}{c|}{Deep Nitsche} & PFNN \\
    \hline
     Unknowns & \multicolumn{2}{c|}{811} & \multicolumn{2}{c|}{811} & 742 \\
    \hline
                    & $\beta=100$ & 0.454\%$\pm$0.072\% & $\beta=100$ & 0.535\%$\pm$0.052\% & \textbf{0.288\%$\pm$0.030\%} \\
    $L=5$           & $\beta=300$ & 1.763\%$\pm$0.675\% & $\beta=300$ & 1.164\%$\pm$0.228\% \\
                    & $\beta=500$ & 5.245\%$\pm$1.943\% & $\beta=500$ & 3.092\%$\pm$1.256\% \\
    \hline
                    & $\beta=100$ & 0.747\%$\pm$0.101\% & $\beta=100$ & 0.483\%$\pm$0.095\% & \textbf{0.309\%$\pm$0.064\%} \\
    $L=6$           & $\beta=300$ & 3.368\%$\pm$0.690\% & $\beta=300$ & 0.784\%$\pm$0.167\% \\
                    & $\beta=500$ & 4.027\%$\pm$1.346\% & $\beta=500$ & 2.387\%$\pm$0.480\% \\
    \hline
                    & $\beta=100$ & 0.788\%$\pm$0.041\% & $\beta=100$ & 0.667\%$\pm$0.149\% & \textbf{0.313\%$\pm$0.071\%} \\
    $L=7$           & $\beta=300$ & 2.716\%$\pm$0.489\% & $\beta=300$ & 1.527\%$\pm$0.435\% \\
                    & $\beta=500$ & 4.652\%$\pm$1.624\% & $\beta=500$ & 1.875\%$\pm$0.653\% \\
    \tabucline[1pt]{-}
  \end{tabu}
\end{table}

In the experiment, we set the essential boundary $\Gamma_D$ to be the left half part of the boundary
and gradually increase the fractal level $L$ of the Koch polygon. With $L$ increased,
the number of segments on the boundary of the domain grows exponentially \cite{koch1904courbe},
posing severe challenges to classical approaches such as the finite element method.
Thanks to the meshfree nature, neural network methods could be more suitable to tackle this problem.
We investigate the performance of the Deep Ritz, Deep Nitsche and PFNN methods with $L=5$, $6$ and $7$.
The configurations of the neural networks are the same to those in the previous experiment.
We report both mean values and standard deviations of the solution errors of all the approaches in Table \ref{T_BVP_2D_MinSurface_Snow},
from which we can see that PFNN can achieve the most accurate results with the least number of parameters
and is much less susceptible to the change of domain boundary.

\subsection{$p$-Liouville-Bratu equation on the Stanford Bunny}

The next test case is a Dirichlet boundary-value problem governed by a $p$-Liouville-Bratu equation:
\begin{equation}
  \label{PDE_p_Bratu}
  -\nabla \cdot \left(|\nabla u|^{p-2} \nabla u \right) - \lambda\exp(u) + c = 0,
\end{equation}
where $p>1$ and $\lambda\geq 0$.
This equation can be transformed to a minimization problem of the following energy functional:
\begin{equation}
  \label{EF_p_Bratu}
  I[w]:= \int_{\Omega}
  \left( \frac{1}{p} |\nabla w|^p - \lambda\exp(w) + cw \right) d\bm{x}.
\end{equation}
The computational domain is a famous 3D graph - the Stanford Bunny \cite{StanfordBunny}.
In this case, the boundary is not given explicitly but determined approximately by 69,451 triangular faces formed by 35,947 vertices.
The shape of the Stanford Bunny is shown in Figure \ref{fig:Test_case_domain} (b),  which is enlarged and translated for better view.

When $p=2$, the $p$-Liouville-Bratu equation is degenerated to the well-known
Liouville-Bratu equation \cite{bratu1914equations}.
The nonlinearity of the problem is increased as $p$ deviates from $2$.
In particular, if $1<p<2$, the diffusivity $|\nabla u|^{p-2}\rightarrow \infty$ as $|\nabla u|\rightarrow 0$.
And if $p>2$ (especially for $p>3$), the diffusivity $|\nabla u|^{p-2}$ grows dramatically as $|\nabla u|$ increases.
In both cases, the corresponding $p$-Liouville-Bratu equation is difficult to solve,
let alone that the computational domain is also very complex.
In addition, the nonlinearity of the problem is also raised as the Bratu parameter $\lambda$ becomes larger.
One thing worth noting is that here $h(u)=-\lambda\exp(u)+c$ is a decreasing function.
Arbitrarily large $\lambda$ could make the solution of (\ref{PDE_p_Bratu}) not be the minimal point
of the energy functional (\ref{EF_p_Bratu}). However, if $\lambda$ is restricted in a proper range,
it is still feasible to obtain a solution by minimizing the equivalent optimization problem.

In the experiment, we set the exact solution to be $u = ((x_1^3-3x_1)(x_2^3-3x_2)-3) (x_3^3+3)$
and conduct two groups of tests to investigate the performance of various neural network methods
in solving the $p$-Liouville-Bratu equation (\ref{PDE_p_Bratu}) with parameter $p=1.2$ and $p=4.0$, respectively.
In each group of tests, we examine the influence of the Bratu parameter with
$\lambda=0.6$ and $\lambda=1.2$, respectively.
Again, the configurations of the neural networks are the same to those in the previous experiment.
The test results are reported in Table \ref{T_BVP_3D_p_Bratu_Bunny}.
From the table, we can clearly see that the proposed PFNN method can outperform the other two methods with less parameters used. 
In particular, for the case of $p=4.0$, both Deep Ritz and Deep Nitsche suffer greatly from the high nonlinearity of the problem and could not produce accurate results,
while PFNN can perform equally well with the change of both $p$ and $\lambda$.

\begin{table}[!h]
  \caption{Errors of various methods for solving the $p$-Liouville-Bratu equation on the Stanford Bunny.}
  \label{T_BVP_3D_p_Bratu_Bunny}
  \centering
  \small
  \renewcommand{\arraystretch}{1.0}
  \begin{tabu}{c|c|l l|l l|c}
    \tabucline[1pt]{-}
    \multicolumn{2}{c|}{Method}   & \multicolumn{2}{c|}{Deep Ritz} & \multicolumn{2}{c|}{Deep Nitsche} & PFNN \\
    \hline
    \multicolumn{2}{c|}{Unknowns} & \multicolumn{2}{c|}{821} & \multicolumn{2}{c|}{821} & 762 \\
    \hline
    \multirow{6}*{$p=1.2$}
     &               & $\beta=100$ & 0.612\%$\pm$0.213\% & $\beta=100$ & 0.659\%$\pm$0.129\% & \textbf{0.513\%$\pm$0.116\%} \\
     & $\lambda=0.6$ & $\beta=300$ & 0.540\%$\pm$0.153\% & $\beta=300$ & 0.593\%$\pm$0.136\% \\
     &               & $\beta=500$ & 0.560\%$\pm$0.218\% & $\beta=500$ & 0.563\%$\pm$0.122\% \\
    \cline{2-7}
     &               & $\beta=100$ & 0.555\%$\pm$0.120\% & $\beta=100$ & 0.643\%$\pm$0.135\% & \textbf{0.489\%$\pm$0.121\%} \\
     & $\lambda=1.2$ & $\beta=300$ & 0.513\%$\pm$0.085\% & $\beta=300$ & 0.608\%$\pm$0.109\% \\
     &               & $\beta=500$ & 0.532\%$\pm$0.159\% & $\beta=500$ & 0.584\%$\pm$0.098\% \\
    \hline
    \multirow{6}*{$p=4.0$}
     &               & $\beta=100$ & 27.646\%$\pm$0.310\% & $\beta=100$ & 28.548\%$\pm$2.849\% & \textbf{0.699\%$\pm$0.467\%} \\
     & $\lambda=0.6$ & $\beta=300$ & 16.327\%$\pm$0.294\% & $\beta=300$ & 21.236\%$\pm$1.326\% \\
     &               & $\beta=500$ & 12.034\%$\pm$0.538\% & $\beta=500$ & 17.972\%$\pm$2.020\% \\
    \cline{2-7}
     &               & $\beta=100$ & 25.133\%$\pm$0.823\% & $\beta=100$ & 30.375\%$\pm$2.387\% & \textbf{0.722\%$\pm$0.393\%} \\
     & $\lambda=1.2$ & $\beta=300$ & 16.330\%$\pm$0.484\% & $\beta=300$ & 21.938\%$\pm$2.369\% \\
     &               & $\beta=500$ & 11.573\%$\pm$0.458\% & $\beta=500$ & 18.998\%$\pm$1.872\% \\
    \tabucline[1pt]{-}
  \end{tabu}
\end{table}

\subsection{Poisson-like equation on a 100D hypercube}
In the last experiment, consider a mixed boundary-value problem
governed by a Poisson-like equation ($p>1$, $\alpha=2$):
\begin{equation}
  \label{BVP_p_Helmholtz}
  \left \{
    \begin{array}{r l}
      -\nabla \cdot \Big(|\nabla u|^{p-2} \nabla u \Big) + u + c = 0,
      & \mbox{in}\ \Omega, \\[1mm]
      u = \varphi,
      & \mbox{on}\ \Gamma_D, \\[1mm]
      \Big( |\nabla u|^{p-2}\nabla u \Big)\cdot \bm{n} + \alpha u = \psi,
      & \mbox{on}\ \Gamma_R,
    \end{array}
  \right.
\end{equation}
where $\Omega$ is a 100D hypercube $[0,1]^{100}$ and $\Gamma_R$ represents the Robin boundary.
For this problem, the corresponding energy functional to minimize is:
\begin{equation}
  \label{EF_p_Helmholtz}
  I[w]:= \int_{\Omega}
  \left( \frac{1}{p} |\nabla w|^p + \frac{1}{2} w^2 + cw \right) d\bm{x}
  + \int_{\Gamma_R} \left( \frac{1}{2}\alpha w^2 - \psi w \right) d\bm{x}.
\end{equation}

In the experiment we set the exact solution to
$$
u = \frac{1}{100}
  \left( \sum_{j=1}^{99} \exp\left(-(x_j^2+x_{j+1}^2)\right) +
         \exp\left(-(x_{100}^2+x_1^2)\right) \right)
$$
and the Dirichlet boundary to be
$\Gamma_D = \cup_{j=1}^{100} \{\bm{x}| x_j=0\}$.
Three groups of tests are carried out, with parameter $p=1.2$, $3.6$ and $4.8$, respectively.
For $p=1.2$ and $p=3.6$, the maximum number of iterations is set to 20,000 epochs,
while for the case that $p=4.8$, it is increased to $40,000$ epochs
since the problem becomes more difficult.
At each step of the iteration, $5,000$ points in the domain and $50$ points
on each hyper-plane that belongs to the boundary are sampled to form the training set.

\begin{figure}[!h]
  \centering
  \subfigure[$p=1.2$]
  {\includegraphics[width=0.32\textwidth]{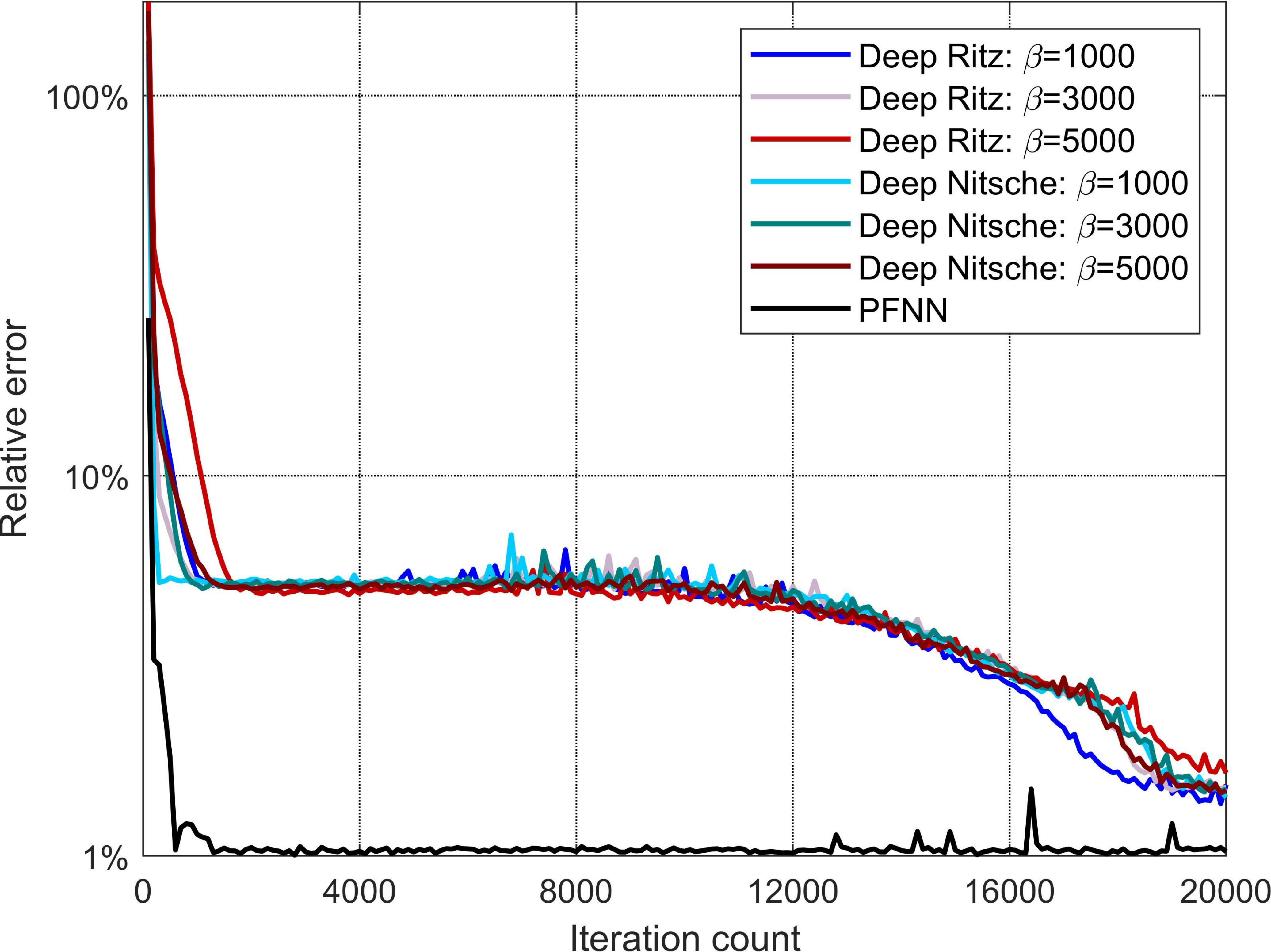}}
  \subfigure[$p=3.6$]
  {\includegraphics[width=0.32\textwidth]{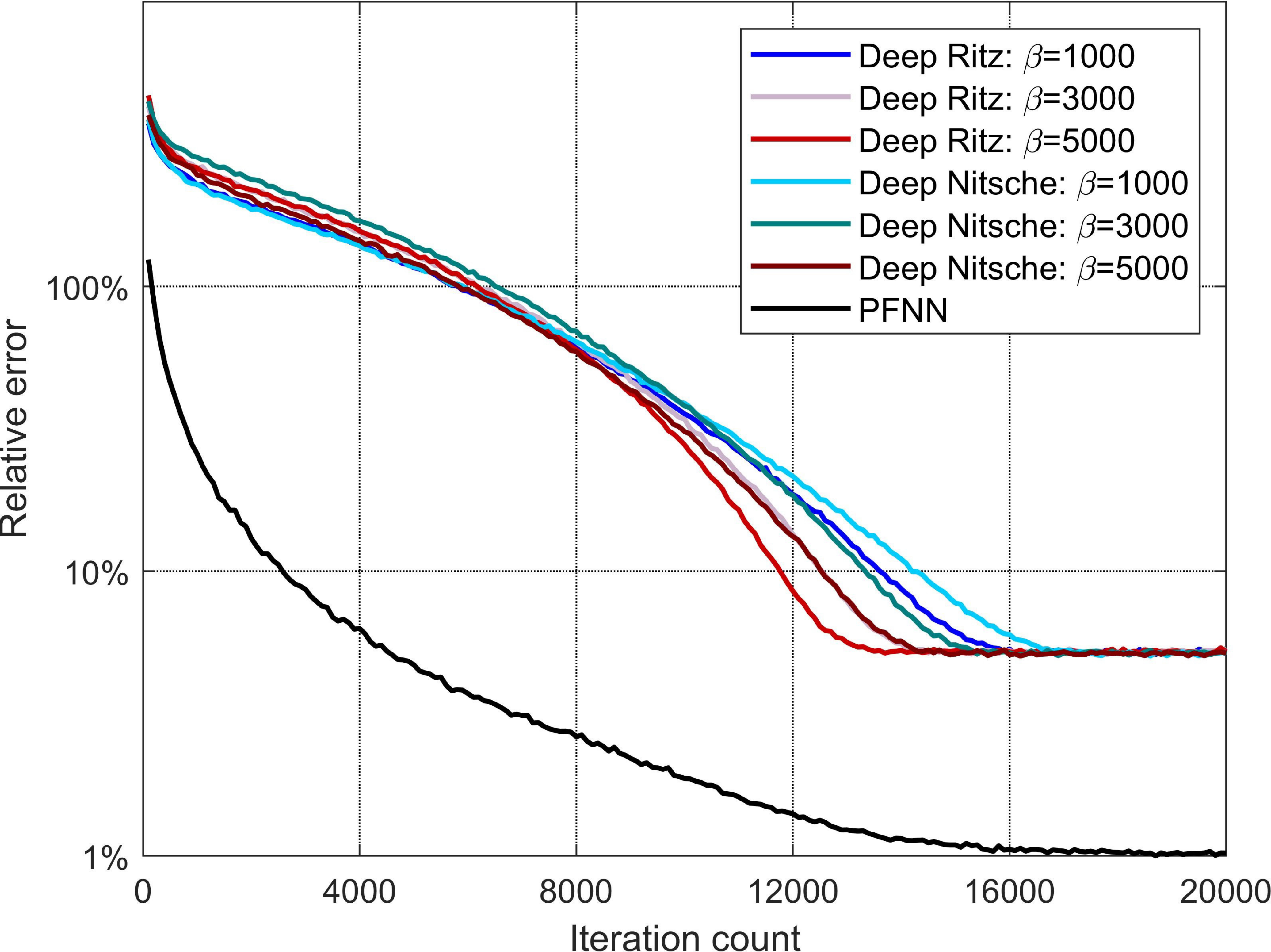}}
  \subfigure[$p=4.8$]
  {\includegraphics[width=0.32\textwidth]{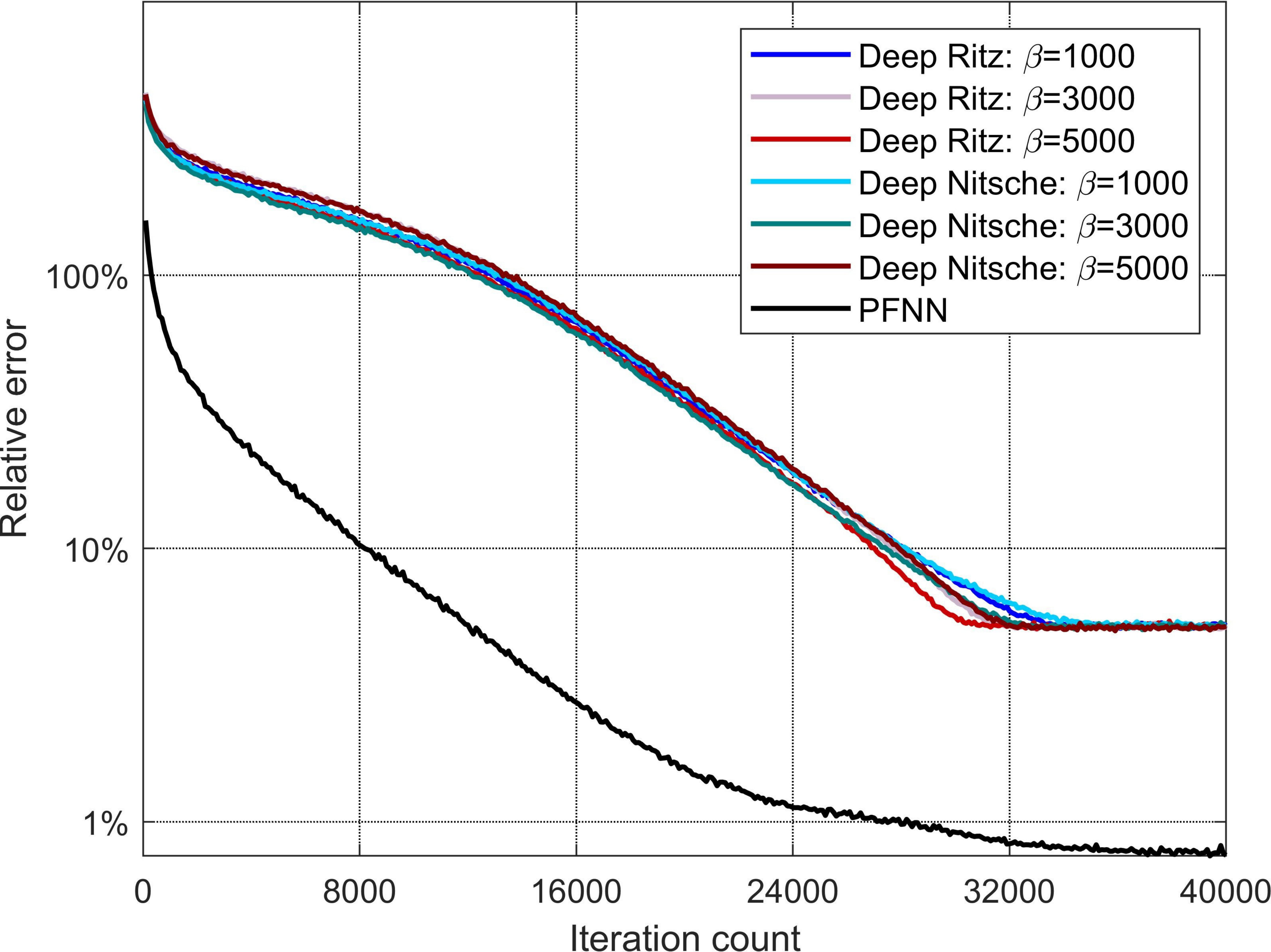}}
  \caption{Error evolution history of various methods during the training process.}
  \label{F_BVP_100D_p_Helmholtz_Errors}
\end{figure}

We apply the Deep Ritz, Deep Nitsche and PFNN methods and compare their performance on this high-dimensional problem.
The network structures for the Deep Ritz and Deep Nitsche methods are the same, both consisting
of $4$ ResNet blocks of width $120$, resulting in $113,881$ undecided parameters.
And the PFNN method employs two neural networks that are comprised of
$1$ and $3$ blocks of width $120$, respectively, corresponding to a total of $26,761+84,841=111,602$ unknowns.
The penalty factors in Deep Ritz and Deep Nitsche are set to three typical values:
$\beta = 1000$, $3000$ and $5000$.
The relative errors of the tree tested approaches are listed in Table \ref{T_BVP_100D_p_Helmholtz},
which again demonstrate the advantages of the proposed PFNN method in terms of both accuracy and robustness.
To further examine the efficiency of various methods, we draw the evolution history of the relative
errors in Figure \ref{F_BVP_100D_p_Helmholtz_Errors}.
The figure clearly indicates that PFNN is not only more accurate, but also has
a faster convergence speed across all tests.

\begin{table}[!h]
  \caption{Errors of various methods for solving the Poisson-like equation on a 100D hypercube.}
  \label{T_BVP_100D_p_Helmholtz}
  \centering
  \small
  \renewcommand{\arraystretch}{1.0}
  \begin{tabu}{c|l l|l l|c}
    \tabucline[1pt]{-}
    Method   & \multicolumn{2}{c|}{Deep Ritz} & \multicolumn{2}{c|}{Deep Nitsche} & PFNN \\
    \hline
    Unknowns & \multicolumn{2}{c|}{113,881} & \multicolumn{2}{c|}{113,881} & 111,602 \\
    \hline
        &          $\beta=1,000$ & 1.405\%$\pm$0.091\% & $\beta=1,000$ & 1.424\%$\pm$0.078\% & \textbf{1.026\%$\pm$0.083\%} \\
    $p=1.2$ &          $\beta=3,000$ & 1.529\%$\pm$0.088\% & $\beta=3,000$ & 1.498\%$\pm$0.082\% \\
        &          $\beta=5,000$ & 1.655\%$\pm$0.084\% & $\beta=5,000$ & 1.487\%$\pm$0.086\% \\
    \hline
        &          $\beta=1,000$ & 5.192\%$\pm$0.062\% & $\beta=1,000$ & 5.199\%$\pm$0.071\% & \textbf{1.037\%$\pm$0.072\%} \\
    $p=3.6$ & $\beta=3,000$ & 5.266\%$\pm$0.055\% & $\beta=3,000$ & 5.133\%$\pm$0.063\% \\
        &          $\beta=5,000$ & 5.203\%$\pm$0.058\% & $\beta=5,000$ & 5.165\%$\pm$0.069\% \\
    \hline
        &          $\beta=1,000$ & 5.206\%$\pm$0.049\% & $\beta=1,000$ & 5.263\%$\pm$0.057\% & \textbf{0.780\%$\pm$0.067\%} \\
    $p=4.8$ &          $\beta=3,000$ & 5.180\%$\pm$0.061\% & $\beta=3,000$ & 5.169\%$\pm$0.044\% \\
        &          $\beta=5,000$ & 5.214\%$\pm$0.053\% & $\beta=5,000$ & 5.187\%$\pm$0.055\% \\
    \tabucline[1pt]{-}
  \end{tabu}
\end{table}

\section{Conclusion}
In this paper, we proposed PFNN -- a penalty-free neural network method to solve a class of
second-order boundary-value problems on complex geometries.
By using the PFNN method the original problem is transformed to a weak form without any penalty terms
and the solution is constructed with two networks and a length factor function.
We provide a theoretical analysis to prove that PFNN is convergent as
the number of hidden units of the networks increases and present
a series of numerical experiments to demonstrate that the PFNN method
is not only more accurate and flexible, but also more robust than the previous state-of-the-art.
Possible future works on PFNN could include further theoretical analysis of the convergence speed and stability, further improvement on the training technique and further study on the corresponding parallel algorithm for solving larger problems on high-performance computers.

\section*{Acknowledgements}
This study was funded in part by Natural Science Foundation of Beijing Municipality (\#JQ18001),
Key-Area R\&D Program of Guangdong Province (\#2019B121204008), and Beijing Academy of Artificial Intelligence.

\bibliographystyle{model1-num-names}
\bibliography{refs}

\end{document}